\documentclass[a4paper]{article}
\usepackage{amsfonts}
\usepackage{amsmath}
\usepackage{amssymb}
\begin{document}
\newtheorem{theorem}{Theorem}[section]
\newtheorem{lemma}[theorem]{Lemma}
\newtheorem{corollary}[theorem]{Corollary}
\newtheorem{conjecture}[theorem]{Conjecture}
\newtheorem{remark}[theorem]{Remark}
\newtheorem{definition}[theorem]{Definition}
\newtheorem{problem}[theorem]{Problem}
\newtheorem{example}[theorem]{Example}
\newtheorem{proposition}[theorem]{Proposition}
\title{{\bf Extension of log pluricanonical forms from subvarieties
}}
\date{October 30, 2007}
\author{Hajime TSUJI}
\maketitle
\begin{abstract}
\noindent In this paper, I prove a very general extension theorem for log pluricanonical  systems.  The strategy and the techniques used here are the same as those in \cite{tu3,tu6,tu7,tu8}. 
The main application of this extension theorem is (together with  Kawamata's subadjunction theorem (\cite{ka3})) to give an optimal subadjunction theorem  which relates the positivities of canonical bundle of the ambient 
projective manifold and that of the (maximal) center of log canonical singularities.  This is an extension of the corresponding result in \cite{tu7}, where I 
dealt with  log pluricanonical systems of general type.    
This subadjunction theorem indicates an approach to solve the abundance conjecture for canonical divisors (or log canonical divisors) in terms of the induction in dimension. 
 \vspace{3mm}  \\
2000 Mathematics Subject Classification: 14J40, 32J18, 32H50
\end{abstract}
\tableofcontents
\section{Introduction}
\noindent In this paper, I present a proof of  the extension theorem of log pluricanonical forms announced  in \cite{tu4,tu6}. The special case of the extension theorem 
has already been proven and used in \cite{tu5,tu6}(cf. \cite[Theorems 2.24,2.25]{tu6}). 
  Although the scheme of the proof is very similar to that of \cite{tu3}, 
it requires a lot more estimates of Bergman kernels and technicalities.

\subsection{Abundance conjecture}

The main motivation to prove such an extension theorem is to investigate 
the pluri (log) canonical systems on a projective varieties.
Since the finite generation of canonical rings has  
been settled very recently (\cite{b-c-h-m}), the most outstanding conjecture in this direction is the following conjecuture. 

\begin{conjecture}{\em ({\bf Abundance conjecture})}\label{abundance} 
Let $X$ be a smooth projective variety defined over $\mathbb{C}$. 
Then $K_{X}$ is abundant, i.e., 
\[
\mbox{\em Kod} (X) = \nu (X)
\]
holds, where $\mbox{\em Kod} (X)$ denotes the Kodaira dimension of $X$ and 
$\nu (X)$ denotes the numerical Kodaira dimension of $X$ (cf. Definition \ref{L-dim}).  $\square$ 
\end{conjecture}
Let us explain the geometric meaning of Conjectre \ref{abundance}.
If $\nu (X) = - \infty$, it is clear that $\mbox{Kod}(X) = -\infty$, 
since $\mbox{Kod}(X) \leqq \nu (X)$ (cf. Definition \ref{L-dim}) always holds. 
Hence in this case, Conjecture \ref{abundance} is trivial.  
Next let $X$ be a smooth projective variety with $\nu (X) \geqq 0$.  
Then $K_{X}$ is pseudoeffective.   
Suppose that $\mbox{Kod}(X) = \nu (X)$ holds. In this case for a sufficiently large
$m$, the rational map associated with $\mid\!\!m!K_{X}\!\!\!\mid$ gives a 
rational fibration (called the Iitaka fibration)
\[
\Phi : = \Phi_{\mid m!K_{X}\!\mid}: X -\cdots\rightarrow Y \subseteq \mathbb{P}^{N_{m}} \hspace{3mm} (N_{m} := \dim \mid\!m!K_{X}\!\!\mid)
\]
with $\dim Y = \nu (X)$. 
Then a general fiber $F$  of $\Phi$ is a smooth projective variety such that 
\[
\mbox{Kod}(F) = \nu (F) = 0.
\]
Let $h$ be an AZD of minimal singularities on $K_{X}$ (cf. Definition \ref{minAZD}).  $(K_{X},h)$ is considered to be the maximal positive part of $K_{X}$. 
Then we see that the curvature current $\Theta_{h}$ of $h$ has no absolutely continuous part on $F$ or equivalently $(K_{X},h)$ is numerically trivial on $F$ 
and $(K_{X},h)\cdot C$ is strictly positive for every irreducible curve $C \subset X$ such that 
$\Phi_{*}(C)\neq 0$, i.e., $\Phi$ is a numerically trivial fibration of $(K_{X},h)$.   
This means that the fibration $\Phi : X -\cdots\rightarrow Y$ extracts all 
the positivity of $(K_{X},h)$. 
In this sense, Conjecture \ref{abundance} asserts that the Iitaka fibration extracts all the positivity of $K_{X}$.  This is the meaning of Conjecture \ref{abundance}.  \vspace{3mm} \\

So far Conjecture \ref{abundance} has been proven only in the case of $\dim X \leqq 3$. 
The main reason why we cannot proceed beyond  the case of $\dim X = 3$ is that 
we do not know how to construct sections of pluricanonical bundles 
without bigness.  
Actually the proof of the abundance conjecture in the case of surfaces depends on
the classification of projective surfaces. 
And the proof of Conjecture \ref{abundance} in the case of projective 
3-folds (\cite{m2,m3,ka2,k-m-m}), the key ingredient of the proof is the clever use of Miyaoka-Yau
type inequality (\cite{m1}).  
But unfortunately this method works well only for the case of 3-folds,
 because the method depends on the speciality of the Riemann-Roch theorem 
 in dimension 3.   
So these partial affirmative solutions of Conjecture \ref{abundance}
do not lead us any further.
Hence it will be desirable to find a systematic method to study pluricanonical systems which works in all dimensions.

Let us consider what is needed to solve Conjecture \ref{abundance}.
For simplicity, first we shall assume that $X$ is already minimal.
Then the Conjecture \ref{abundance} is equivalent to the stable base point 
freeness of $K_{X}$ by a theorem of Kawamata (\cite{freeness}). 
Thus the abundance conjecture can be viewed as a base point freeness theorem.

To prove the base point freeness of some linear systems, especially 
the case of adjoint line bundles,   the well known method is 
to produce log canonical center and go by induction in dimension
using Kawamata-Viehweg vanishing theorem. 

To solve Conjecture \ref{abundance}, we shall consider a similar approach. 
Let $X$ be a smooth projective variety with pseudoeffective canonical bundle
and let $h$ be an AZD of $K_{X}$ with minimal singularities.  
The strategy is as follows. 
\begin{enumerate}
\item Find a LC center (log canonical center), say  $S$ for $\alpha K_{X}$ for some 
$\alpha > 0$.  Here the LC center means a little bit broader sense, i.e., 
LC center here means the cosupport of the multiplier ideal sheaf with 
respect to a singular hermitian metric of $\alpha K_{X}$ with  
semipositive curvature current.  We shall assume that $S$ is smooth and is a maximal LC center (This can be assured by taking an embedded resolution
of the LC center).   
\item Use subadjunction theorem due to Kawamata (\cite{ka3})  
to compare the canonical divisor of the center and $K_{X}$.  
\item  Lift  the pluricanonical system of the center (possibly twisted by some fixed ample line bundle $B$ on $X$)   to 
a  sub linear system of the pluricanonical system of $X$ (possibly twisted by 
the ample line bundle $B$). 
\end{enumerate}  

\noindent 
If this strategy works, we can see  the positivity of $(K_{X},h)\!\!\mid_{S}$ dominates that of $K_{S}$.   If we have already known  the abundance of $K_{S}$ and 
the equality $\nu (S) = \nu (K_{X}\!\!\mid_{S},h\!\!\mid_{S})$ holds,   
then we see that the positivity of  $(K_{X},h)\!\!\mid_{S}$ can be extracted 
in terms of the Iitaka fibration of $S$. 
Moreover if we are able to construct the LC center as above through every point of $X$ (this is expected when $\nu (X) \geqq 1$), 
then we may see that for a general such $S$ above have pseudoeffective $K_{S}$
and the equality $\nu (S) = \nu (K_{X}\!\!\mid_{S},h\!\!\mid_{S})$ holds.  
As a consequence, we may extract the positivity of $(K_{X},h)$ in terms of 
a family of LC centers.  And it is not difficult to prove Conjecture \ref{abundance} in this case, if we have already proven  Conjecture \ref{abundance} 
for every smooth projective variety of dimension $< \dim X$. 
Hence this strategy can be viewed as an approach 
to Conjecture \ref{abundance} by the induction in dimension. 
\vspace{3mm} \\ 
 
This strategy has been first considered by the author in the series of papers \cite{tu5,tu6,tu7}.  But in these papers, the varieties are assumed to be  of general type and in this case the LC center $S$ is also of general type.  
Hence what we compare is the volume of $K_{S}$ and that of $(K_{X},h)\!\!\mid_{S}$ and  there are essentially no analytic difficulties.   In fact \cite{h-m,ta}
have interpreted \cite{tu7} in terms of algebro geometric languages. 

But in contrast to the case of general type, the case of non general type 
is much harder in analysis.  The reason is that there is no room of 
positivity to approximate AZD of $K_{X}$ by a sequence of pseudoeffective 
singular hermitian metrics with algebraic singularities. 
This phenomena (loss of positivity) was first observed in the work of 
Demailly (\cite{dem}) on the regularization of closed positive currents.

Hence in the case of non general type, we need to provide effective 
estimates to obtain the desired singular hermitian metrics  
as was observed in \cite{tu3,si}. 
In this case the nefness should be replaced by the {\bf semipositivity 
of curvature}.   
This is the crucial point to apply the $L^{2}$-extension theorem 
which is considered to be the substitute of Kawamata-Viehweg vanishing 
theorem (\cite{ka1}) in the case of varieties of non general type. 
We summerize the comparison of the cases of general type and of non general type in the following table.\vspace{5mm} \\
\hspace{-10mm}
\begin{tabular}{|l||c|c|}
\hline 
 & \textbf{Case of general type} & \textbf{Case of non general type} \\
 \hline
 To produce LC center & concentration of multiplicities & ?   \\
 \hline
 Necessary positivity &  nefness & curvature semipositivity \\
 \hline 
 To lift sections & Kawamata-Viehweg vanishing theorem & $L^{2}$-extension theorem    \\     
\hline 
\end{tabular}
\vspace{5mm}\\
The purpose of this article is to implement  the steps 2 and 3 of the strategy above.   
At this moment I do not know how to produce LC centers.    
\subsection{Main results} 

Let us explain the main results in this article. 
Let $X$ be a smooth projective variety 
and let $(L,h_{L})$ be a singular hermitian line bundle on $X$ such that 
$\Theta_{h_{L}}\geqq 0$\footnote{Here we have used the convention 
that $\Theta_{h_{L}} = \sqrt{-1}\bar{\partial}\partial\log h_{L}$.
The advantage of this convention is that $\Theta_{h_{L}}$ is always a real current.} on $X$.  
We assume that $h_{L}$ is lowersemicontinuous. 
Throughout this paper we shall assume that all the singular hermitian metrics 
are lowersemicontinuous.  

Let $m_{0}$ be a positive integer. 
Let $\sigma_{0} \in \Gamma (X,{\cal O}_{X}(m_{0}L)\otimes {\cal I}_{\infty}(h_{L}^{m_{0}}))$ be a  bounded 
global section (cf. Section \ref{singhm} for the definition of 
${\cal I}_{\infty}(h_{L}^{m_{0}})$). 
Let $\alpha$ be a positive rational number $\leqq 1$ and let $S$ be 
an irreducible component of the center of LC(log canonical) singularity  but not KLT(Kawamata log terminal) and $(X,(\alpha -\epsilon )(\sigma_{0}))$ is KLT on the generic point of $S$ 
for every $0 < \epsilon << 1$. 
We set 
\[
\Psi_{S} = \alpha\cdot\log h_{L}^{m_{0}}(\sigma_{0} ,\sigma_{0}).
\]
Suppose that $S$ is smooth for simplicity.  
Let $dV$ be a $C^{\infty}$ volume form on $X$. 

In this situation we may define a (possibly singular) measure 
$dV[\Psi_{S} ]$ on $S$ as the residue as follows.   
Let $f : Y \longrightarrow X$ be a log resolution of 
$(X,\alpha (\sigma_{0} ))$. 
Then we may define the singular volume form $f^{*}dV[f^{*}\Psi_{S} ]$ 
on the divisorial components of $f^{-1}(S)$ (the volume form is identically 
$0$ on the components with discrepancy $> -1$) by taking residue 
along $f^{-1}(S)$. 
The singular volume form $dV[\Psi_{S} ]$ is defined as the fibre integral of 
$f^{*}dV[f^{*}\Psi_{S} ]$ (the actual integration takes place only on the components with discrepancy $-1$), i.e., $dV[\Psi_{S}]$ is the  
residue volume form of general codimension. 

Let $dV_{S}$ be a $C^{\infty}$ volume form on $S$ and 
let $\varphi$ be the function on $S$ defined by
\[
\varphi := \log \frac{dV_{S}}{dV[\Psi_{S} ]}
\]
($dV[\Psi_{S} ]$ may be singular on a subvariety of $S$, also 
it may be totally singular on $S$).   
The following is the main theorem in this article.

\begin{theorem}\label{main}(cf. \cite[Theorem 5.1]{tu5})
Let $X$,$S$,$\Psi_{S}$ be as above.  
Suppose that $S$ is smooth. 
Let $d$ be a positive integer such that $d > \alpha m_{0}$. 
We assume  that $(K_{X}+dL,e^{-\varphi}\cdot (dV^{-1}\cdot h_{L}^{d})\mid_{S})$
 is weakly pseudoeffective (cf. Definition \ref{wpe}) and  let  $h_{S}$ be an AZD of $(K_{X}+dL,e^{-\varphi}\cdot (dV^{-1}\cdot h_{L}^{d})\mid_{S})$.
  
Then every element of 
\[
H^{0}(S,{\cal O}_{S}(m(K_{X}+dL))\otimes{\cal I}(e^{-\varphi}\cdot h_{L}^{d}\mid_{S}\cdot h_{S}^{m-1}))
\]
extends to an element of 
\[
H^{0}(X,{\cal O}_{X}(m(K_{X}+dL))\otimes{\cal I}(h_{L}^{d})). 
\]
In particular every element of 
\[
H^{0}(S,{\cal O}_{S}(m(K_{X}+dL))\otimes{\cal I}(e^{-\varphi}\cdot h_{L}^{d}\mid_{S})\cdot {\cal I}_{\infty}(e^{-(m-1)\varphi}\cdot h^{m-1}_{L}\mid_{S}))
\] 
extends to an element of 
\[
H^{0}(X,{\cal O}_{X}(m(K_{X}+dL))\otimes{\cal I}(h_{L}^{d}\cdot h_{0}^{m-1})). 
\]
where $h$ is an AZD of $K_{X} + dL$ of  minimal singularities. 
$\square$ \end{theorem}

\noindent As we mentioned as above the smoothness assumption on $S$ is 
just to make the statement simpler.   
And it may be worthwhile to note that the weight function $\varphi$ is not 
necessary when $dV[\Psi_{S}]$ is locally $L^{1}$ on $S$ and $h_{L}$ is 
bounded on $S$ (see the proof in Section 7). 

Theorem \ref{main} follows from  Theorem \ref{subad2}
below by using a limiting process (cf. Section 9).

\begin{theorem}\label{subad2}
Let $X$,$S$,$\Psi_{S}$ be as above.  
Suppose that $S$ is smooth.
Let $(E,h_{E})$ be a pseudoeffective singular hermitian line bundle on $X$
(cf. Definition \ref{pe}).  
Let $d$ be a positive integer such that $d > \alpha m_{0}$.
Let $m$ be a positive integer. 
We assume that \\
$(K_{X}+dL+\frac{1}{m}E\mid_{S},e^{-\varphi}\cdot (dV^{-1}\cdot h_{L}
\cdot h_{E}^{\frac{1}{m}})\mid_{S})$
is weakly pseudoeffective.   
Let $h_{S}$ be an AZD of $(K_{X}+dL+\frac{1}{m}E\mid_{S},e^{-\varphi}\cdot (dV^{-1}\cdot h_{L}^{d}\cdot h_{E}^{\frac{1}{m}})\mid_{S})$. 
Suppose that $h_{S}$ is normal (cf. Definition \ref{normal}) or
\[
\dim H^{0}(S,{\cal O}_{S}(A+m\ell(K_{X}+dL+\frac{1}{m}E))\mid_{S})\otimes{\cal I}(h_{S}^{m\ell})) = O(\ell^{\nu})
\]
holds for every ample line bundle $A$ on $X$, where $\nu$ denotes the numerical Kodaira dimension 
$\nu_{num}(K_{X}+L+\frac{1}{m}E\mid_{S},h_{S})$ 
of $(K_{X}+L+\frac{1}{m}E\mid_{S},h_{S})$ \\(cf. Definition \ref{numerical}).
\\\\\
Then every element of 
\[
H^{0}(S,{\cal O}_{S}((m+1)(K_{X}+dL)+E)\otimes {\cal I}(e^{-\varphi}\cdot h_{L}^{d}\mid_{S}\cdot  h_{S}^{m}))
\]
extends to an element of 
\[
H^{0}(X,{\cal O}_{X}((m+1)(K_{X}+dL)+E)\otimes {I}(h_{L}^{d}\cdot h_{0,\frac{1}{m}}^{m})), 
\]
where $h_{0,\frac{1}{m}}$ is an AZD of $K_{X}+dL+\frac{1}{m}E$ of minimal singularities. 
$\square$
\end{theorem}
As an application of Theorem \ref{subad2}, we obtain the following 
theorem. 

\begin{theorem}\label{numsubad}
Let $X$ be a smooth projective variety and let $D$ be an effective 
$\mathbb{Q}$-divisor on $X$ numerically equivalent to $\alpha K_{X}$ 
for some $\alpha > 0$. Let $h$ be the supercanonical AZD  of $K_{X}$ 
(\cite{tu9}) or any AZD of minimal singularities (cf. Definition \ref{minAZD}).   
Let $S$ be a maximal log canonical center of 
$(X,D)$ in the sense that $S$ cannot be a proper subvariety of the log canonical center of $(X,D)$. We assume that for every rational number  $0< \epsilon << 1$,   
$(X,(1-\epsilon )D)$ is log terminal at the generic point of $S$.    
Suppose that $S$ is smooth and $K_{S}$ is pseudoeffective.

Then for every ample line bundle $A$ and every sufficiently ample line bundle B and every positive integer $m$, there exists an injection 
\[
\hspace{-70mm} H^{0}(S,{\cal O}_{S}(mK_{S} + A))\hookrightarrow
\]
\[
\mbox{\em Image}\{ H^{0}(X,{\cal O}_{X}(m(1+d)K_{X}+A + B)\otimes {\cal I}(h^{m-1}))
\rightarrow   
H^{0}(S,{\cal O}_{S}(m(1+d)K_{X} + A+ B))\}.
\]
In particular  
\[
\nu (K_{S}) \leqq \nu_{asym}(K_{X}\!\!\mid_{S},h\!\mid_{S})
\]
holds (for the definitions of $\nu (K_{S})$ and $\nu_{asym}(K_{X}\!\!\!\mid_{S},h\!\!\mid_{S})$
see Definitions \ref{L-dim} and \ref{numerical Kodaira}  below).  
$\square$
\end{theorem}
\begin{remark}
The smoothness assumption of $S$ is not essential in Theorem \ref{numsubad}.  In fact we just need to take an embedded
resolution. 
$\square$ 
\end{remark}
The following theorem is a variant of Theoem \ref{numsubad}. 
\begin{theorem}\label{numsubad2}
Let $X$ be a smooth projective variety and let $D$ be an effective 
$\mathbb{Q}$-divisor on $X$ numerically equivalent to $\alpha K_{X}$ 
for some $\alpha > 0$. Let $h$ be the supercanonical AZD  of $K_{X}$ 
(\cite{tu9}) or any AZD of minimal singularities (cf. Definition \ref{minAZD}).   
Let $S$ be a maximal log canonical center of 
$(X,D)$ in the sense that $S$ cannot be a proper subvariety of the log canonical center of $(X,D)$. We assume that for every rational number  $0< \epsilon << 1$,   
$(X,(1-\epsilon )D)$ is log terminal at the generic point of $S$.    
Suppose that $S$ is smooth and $K_{S}$ is pseudoeffective.

Then for every sufficiently ample line bundle B and every positive integer $m$, there exists an injection 
\[
\hspace{-70mm} H^{0}(S,{\cal O}_{S}(mK_{S}))\hookrightarrow
\]
\[
\mbox{\em Image}\{ H^{0}(X,{\cal O}_{X}(m(1+d)K_{X}+ B)\otimes{\cal I}(h^{m-1}))
\rightarrow   
H^{0}(S,{\cal O}_{S}(m(1+d)K_{X} + B))\}.
\]
$\square$
\end{theorem}
\begin{remark}\label{2rm}
One can deduce Theorem \ref{numsubad2} from Theorem \ref{numsubad} by 
taking $A$ sufficiently ample in Theorem \ref{numsubad}.  
But the proof in Section 9 below implies that we just need to take $B$ 
with a continuous metric with semipositive curvature such  that there exists an 
injection 
\[
H^{0}(S,{\cal O}_{S}(mK_{S})) \hookrightarrow 
H^{0}(S,{\cal O}_{S}(m(1+d)K_{X}+ B)\otimes {\cal I}(h^{m}\!\!\mid_{S}))
\]
exists for all $m\geqq 1$.  
Hence the ampleness of $B$ is somewhat irrelevant in Theorem \ref{numsubad2}.  $\square$ 
\end{remark}  
\noindent Theorem \ref{numsubad} implies that the positivity of 
$K_{S}$ is dominated by the positivity of $K_{X}$ (up to a constant multiple). 
Theorem \ref{numsubad} is  obtained by a similar argument as in \cite{tu5,tu6,tu7} using Theorem \ref{subad1} below and Kawamata's semipositivity theorem (\cite[Theorem 2]{ka3}, see Theorem \ref{pos} below). 

In Theorem \ref{numsubad}, the presence of $h$ is crucial, even if  
$K_{X}$ is assumed to be nef.  Because in this case, the abundace of $K_{X}$ 
implies that we may replace $(K_{X},h)$ by $K_{X}$ in Theorem \ref{numsubad}.
Hence this is somewhat opposite.  

We note that for a minimal algebraic variety $X$, the abundance of $K_{X}$ implies that $K_{X}$ is numerically trivial on the fiber of the Iitaka fibration and 
the pluricanonical system comes from an ample Kawamata log terminal divisor 
on the base of the Iitaka fibration (\cite{f-m}).  
Hence Theorems \ref{numsubad} is a supporting evidence of Conjecture \ref{abundance}.  In fact assuming the existence of minimal models for projective varieties with pseudoeffective canonical bundles and the abundance conjecture, we may 
easily deduce Theorem \ref{numsubad} by using the $L^{2}$-extension theorem (Theorems \ref{o-t} and \ref{extension}). \vspace{3mm}\\

The organization of this article is as follows. 
In Section 2, I collect basic tools to prove Theorem \ref{main}. 
This section is mainly for algebraists.  
In Section 3, I define the numerical Kodaira dimension 
and the asymptotic Kodaira dimension of 
a pseudoeffective singular hermitian line bundle on a smooth projective variety  and study  
the relation between the two dimensions. 
In Section 4,  I relate the Monge-Amp\`{e}re measure of the curvature 
current of a pseudoeffective singular hermitian line bundle to the asymptotic 
expansion of Bergman kernels. 
This leads us to define the local volume of a pseudoeffective singular hermitian line bundles and a natural and very interesting conjecture 
for the relation between the asymptotics of Bergman kernels and the Monge-Amp\`{e}re mass.  This direction should be studied in near future.
In Section 5, I prove  an analogue of  Kodaira's lemma for big pseudoeffective
singular hermitian line bundles.  
This lemma is used in the next section.       
In Section 6, I prove the dynamical construction of an AZD for  adjoint type 
singular hermitian line bundles.  The idea of the proof is similar to  
the one in \cite{tu6}.  But it requires a little bit  more complication. 
In Section 7,  I prove Theorem \ref{subad1} by using the dynamical construction of an AZD in Section 6. 
In Section 8, I prove Theorem \ref{numsubad}   by using 
Kawamata's semipositivity theorem (Theorem \ref{pos}) and Theorem \ref{subad2}.
In Section 9, I prove Theorem \ref{main} combining results in the previous 
sections.

\section{Preliminaries}
In this section we collect the basic tools.  They are standard except 
Definitions \ref{wpe},\ref{normal} and \ref{singAZD}.   
\subsection{Singular hermitian metrics}\label{singh}
In this subsection $L$ will denote a holomorphic line bundle on a complex manifold $X$. 
\begin{definition}\label{singhm}
A  singular hermitian metric $h$ on $L$ is given by
\[
h = e^{-\varphi}\cdot h_{0},
\]
where $h_{0}$ is a $C^{\infty}$ hermitian metric on $L$ and 
$\varphi\in L^{1}_{loc}(X)$ is an arbitrary function on $X$.
We call $\varphi$ a  weight function of $h$. $\square$ 
\end{definition}
The curvature current $\Theta_{h}$ of the singular hermitian line
bundle $(L,h)$ is defined by
\[
\Theta_{h} := \Theta_{h_{0}} + \sqrt{-1}\partial\bar{\partial}\varphi ,
\]
where $\partial\bar{\partial}$ is taken in the sense of a current 
and  we have used the convention
that $\Theta_{h_{0}} = \sqrt{-1}\bar{\partial}\partial\log h_{0}$.
The $L^{2}$ sheaf ${\cal L}^{2}(L,h)$ of the singular hermitian
line bundle $(L,h)$ is defined by
\[
{\cal L}^{2}(L,h)(U) := \{ \sigma\in\Gamma (U,{\cal O}_{X}(L))\mid 
\, h(\sigma ,\sigma )\in L^{1}_{loc}(U)\} ,
\]
where $U$ runs over the  open subsets of $X$.
In this case there exists an ideal sheaf ${\cal I}(h)$ such that
\[
{\cal L}^{2}(L,h) = {\cal O}_{X}(L)\otimes {\cal I}(h)
\]
holds.  We call ${\cal I}(h)$ the {\bf multiplier ideal sheaf} of $(L,h)$.
If we write $h$ as 
\[
h = e^{-\varphi}\cdot h_{0},
\]
where $h_{0}$ is a $C^{\infty}$ hermitian metric on $L$ and 
$\varphi\in L^{1}_{loc}(X)$ is the weight function, we see that
\[
{\cal I}(h) = {\cal L}^{2}({\cal O}_{X},e^{-\varphi})
\]
holds.
For $\varphi\in L^{1}_{loc}(X)$ we define the multiplier ideal sheaf of $\varphi$ by 
\[
{\cal I}(\varphi ) := {\cal L}^{2}({\cal O}_{X},e^{-\varphi}).
\] 
Similarly for $1\leqq p \leqq + \infty$, we define
\[
{\cal L}^{p}(L,h)(U) := \{ \sigma\in\Gamma (U,{\cal O}_{X}(L))\mid 
\, h(\sigma ,\sigma )\in L^{p/2}_{loc}(U)\} ,
\]
where $U$ runs over the  open subsets of $X$.
In this case there exists an ideal sheaf ${\cal I}_{p}(h)$ such that
\[
{\cal L}^{p}(L,h) = {\cal O}_{X}(L)\otimes {\cal I}_{p}(h)
\]
holds.  We call ${\cal I}_{p}(h)$ the {\bf $L^{p}$ multiplier ideal sheaf} of $(L,h)$.

\begin{remark}
It is known that ${\cal I}(h)$ is coherent when $\Theta_{h}$ is locally 
bounded from below by a $C^{\infty}$ form (\cite{n}).
But it is not clear whether  ${\cal I}_{p}(h)$ is coherent under the same
condition for $p\neq 2$. $\square$
\end{remark}

\begin{example}
Let $\sigma\in \Gamma (X,{\cal O}_{X}(L))$ be the global section. 
Then 
\[
h := \frac{1}{\mid\sigma\mid^{2}} = \frac{h_{0}}{h_{0}(\sigma ,\sigma)}
\]
is a singular hemitian metric on $L$, 
where $h_{0}$ is an arbitrary $C^{\infty}$ hermitian metric on $L$
(the right hand side is ovbiously independent of $h_{0}$).
The curvature $\Theta_{h}$ is given by
\[
\Theta_{h} = 2\pi (\sigma ), 
\]
where $(\sigma )$ denotes the current of integration over the 
divisor of $\sigma$. $\square$ 
\end{example}

First we define the pseudoeffectivity of a singular hermitian line bundle. 
\begin{definition}\label{pe}
A line bundle $L$ on a complex manifold is said to be {\bf pseudoeffective}, if there exists 
a singular hermitian metric $h$ on $L$ such that 
the curvature current 
$\Theta_{h}$ is a closed positive current.
A singular hermitian line bundle $(L,h_{L})$ is said to be {\bf pseudoeffective},
if the curvature current $\Theta_{h_{L}}$ is a closed positive current. 
$\square$ 
\end{definition}
An important class of singular hermitian metrics with semipositive
curvature current is  algebraic singular hermitian metrics. 
\begin{definition}\label{algsing}
Let $h$ be a singular hermitian metric on $L$. 
We say that $h$ is {\bf algebraic}, if there exists a positive 
integer $m_{0}$ and global holomorphic sections $\sigma_{0},\cdots ,\sigma_{N}$
of $m_{0}L$ such that 
\[
h= (\sum_{i=0}^{N}\mid\sigma_{i}\mid^{2})^{-\frac{1}{m_{0}}}
\]
holds.   
$\square$
\end{definition}
Definition \ref{algsing} is naturally generalized to the case of 
$\mathbb{Q}$-line bundles in an ovbious way.

Also the following weaker version of pseudoeffectivity is also 
important. 

\begin{definition}\label{wpe}
Let $(L,h_{L})$ be a singular hermitian line bundle on a smooth projective 
variety $X$.  $(L,h_{L})$ is said to be {\bf weakly pseudoeffective}, if there exists an 
ample line bundle $A$ on $X$ such that 
\[
H^{0}(X,{\cal O}_{X}(mL + A)\otimes {\cal I}(h_{L}^{m})) \neq 0
\]
holds for every $m \geqq 0$.
A line bundle on a smooth projective variety is said to be pseudoeffective 
there exists an 
ample line bundle $A$ on $X$ such that 
\[
H^{0}(X,{\cal O}_{X}(mL + A)) \neq 0 
\]
holds for every $m \geqq 0$.
 $\square$
\end{definition}
\begin{remark}
Let $(L,h_{L})$ is a pseudoeffective singular hermitian line bundle on a smooth projective variety $X$.  
Then $(L,h_{L})$ is weakly pseudoeffective.  
This follows from  an easy application of H\"{o}rmander's $L^{2}$-estimate.
 $\square$
\end{remark}

The following definition is useful in this article. 
\begin{definition}\label{normal}
Let $(L,h_{L})$ be a pseudoeffective singular hermitian line bundle on
a smooh projective variety $X$. 
$(L,h_{L})$ is said to be {\bf normal}, if the set  
\[
E := \{ x\in X\mid n(\Theta_{h_{L}},x) > 0\}
\]
is containd in a proper analytic set of $X$, where
$n(\Theta_{h_{L}},x)$ denotes the Lelong number of the closed positive current $\Theta_{h_{L}}$ at 
$x\,\,\,$\footnote{Usually I use $\nu$ instead of $n$.  But in this article, I use $n$ not to confuse with the numerical Kodaira dimension or the asymptotic Kodaira dimension. }. $\square$
\end{definition}
\begin{remark}
By the fundamental theorem of Siu (\cite{lelong}), $E$ is at most 
a countable union of subvarieties of $X$. $\square$ 
\end{remark}

\subsection{Analytic Zariski decompositions (AZD)}\label{subsecAZD}
In this subsection we shall introduce the notion of analytic Zariski decompositions. 
By using analytic Zariski decompositions, we can handle a pseudoeffective line bundles
like  nef line bundles.

\begin{definition}\label{defAZD}
Let $X$ be a compact complex manifold and let $L$ be a holomorphic line bundle
on $X$.  A singular hermitian metric $h$ on $L$ is said to be 
an analytic Zariski decomposition, if the followings hold.
\begin{enumerate}
\item $\Theta_{h}$ is a closed positive current,
\item for every $m\geq 0$, the natural inclusion
\[
H^{0}(X,{\cal O}_{X}(mL)\otimes{\cal I}(h^{m}))\rightarrow
H^{0}(X,{\cal O}_{X}(mL))
\]
is an isomorphim. $\square$
\end{enumerate}
\end{definition}
\begin{remark} If an AZD exists on a line bundle $L$ on a smooth projective
variety $X$, $L$ is pseudoeffective by the condition 1 above. $\square$
\end{remark}

\begin{theorem}(\cite{tu,tu2})
 Let $L$ be a big line  bundle on a smooth projective variety
$X$.  Then $L$ has an AZD. $\square$  
\end{theorem}
As for the existence of AZD for general pseudoeffective line bundles, 
now we have the following theorem.
\begin{theorem}(\cite[Theorem 1.5]{d-p-s})\label{AZD}
Let $X$ be a smooth projective variety and let $L$ be a pseudoeffective 
line bundle on $X$.  Then $L$ has an AZD. $\square$
\end{theorem}
Although the proof is in \cite{d-p-s}, 
we shall give a proof here, because we shall use it afterwards. 

 Let  $h_{0}$ be a fixed $C^{\infty}$ hermitian metric on $L$.
Let $E$ be the set of singular hermitian metric on $L$ defined by
\[
E = \{ h ; h : \mbox{lowersemicontinuous singular hermitian metric on $L$}, 
\]
\[
\hspace{70mm}\Theta_{h}\,
\mbox{is positive}, \frac{h}{h_{0}}\geq 1 \}.
\]
Since $L$ is pseudoeffective, $E$ is nonempty.
We set 
\[
h_{L} = h_{0}\cdot\inf_{h\in E}\frac{h}{h_{0}},
\]
where the infimum is taken pointwise. 
The supremum of a family of plurisubharmonic functions 
uniformly bounded from above is known to be again plurisubharmonic, 
if we modify the supremum on a set of measure $0$(i.e., if we take the uppersemicontinuous envelope) by the following theorem of P. Lelong.

\begin{theorem}(\cite[p.26, Theorem 5]{l})\label{lelong}
Let $\{\varphi_{t}\}_{t\in T}$ be a family of plurisubharmonic functions  
on a domain $\Omega$ 
which is uniformly bounded from above on every compact subset of $\Omega$.
Then $\psi = \sup_{t\in T}\varphi_{t}$ has a minimum 
uppersemicontinuous majorant $\psi^{*}$  which is plurisubharmonic.
We call $\psi^{*}$ the uppersemicontinuous envelope of $\psi$. $\square$ 
\end{theorem}
\begin{remark} In the above theorem the equality 
$\psi = \psi^{*}$ holds outside of a set of measure $0$(cf.\cite[p.29]{l}).
$\square$  
\end{remark}
In this paper, we shall call the uppersemicontinuous envelope (resp. lower semicontinuous envelope)  by the upper envelope (resp. the lower envelope)
for simplicity. 
By Theorem \ref{lelong}, we see that $h_{L}$ is also a 
singular hermitian metric on $L$ with $\Theta_{h}\geq 0$.
Suppose that there exists a nontrivial section 
$\sigma\in \Gamma (X,{\cal O}_{X}(mL))$ for some $m$ (otherwise the 
second condition in Definition 2.3 is empty).
We note that  
\[
\log \mid\sigma\mid^{\frac{2}{m}}
\]
gives the weight of a singular hermitian metric on $L$ with curvature 
$2\pi m^{-1}(\sigma )$, where $(\sigma )$ is the current of integration
along the zero set of $\sigma$. 
By the construction we see that there exists a positive constant 
$c$ such that  
\[
\frac{h_{0}}{\mid\sigma\mid^{\frac{2}{m}}} \geq c\cdot h_{L}
\]
holds. 
Hence
\[
\sigma \in H^{0}(X,{\cal O}_{X}(mL)\otimes{\cal I}_{\infty}(h_{L}^{m}))
\]
holds.  
Hence in particular
\[
\sigma \in H^{0}(X,{\cal O}_{X}(mL)\otimes{\cal I}(h_{L}^{m}))
\]
holds.  
 This means that $h_{L}$ is an AZD of $L$. $\square$  
\begin{remark}\label{linfty}
By the above proof we have that for the AZD $h_{L}$ constructed 
as above,
\[
H^{0}(X,{\cal O}_{X}(mL)\otimes{\cal I}_{\infty}(h_{L}^{m}))
\simeq 
H^{0}(X,{\cal O}_{X}(mL))
\]
holds for every $m$.  $\square$
\end{remark}
It is easy to see that the multiplier ideal sheaves 
of $h_{L}^{m}(m\geq 1)$ constructed in the proof of
 Theorem 2.2 are independent of 
the choice of the $C^{\infty}$ hermitian metric $h_{0}$.
The AZD constructed as in the proof of Theorem \ref{AZD} has minimal singularity in the following sense. 
\begin{definition}\label{minAZD}
Let $L$ be a pseudoeffective line bundle on a smooth projective variety $X$.
An AZD $h$ on $L$  is said to be an {\bf AZD of minimal singularities}, if 
for any AZD $h^{\prime}$ on $L$, there exists a positive constant $C$ such that
\[
h \leqq   C \cdot h^{\prime}
\]
holds. $\square$ 
\end{definition}
\begin{remark} In \cite{tu9}, I have constructed a canonical AZD (the supercanonical AZD) of the canonical bundle of a projective variety with pseudoeffective canonical bundle.  The supercanonical AZD is completely determined by the complex structure.   In the previous papers, I have called an AZD of minimal singularities a canonical AZD.  Since this may cause a confusion, I have changed the name.$\square$ 
\end{remark}
The following proposition is trivial but important. 
\begin{proposition}\label{normality criterion}
Let $h$ be an AZD of a line bundle $L$ on a compact complex manifold $X$. 
Suppose that there exists a positive integer $m_{0}$ such that 
$\mid\!m_{0}L\!\mid\neq\emptyset$, then $h$ is normal. $\square$
\end{proposition}
{\bf Proof.} 
Suppose that there exists a positive integer $m_{0}$ such that 
$\mid\!m_{0}L\!\mid \neq \emptyset$. 
Let $x\in M$ be a point such that $n(\Theta_{h},x) > 0$.  
Then ${\cal I}(h^{m})_{x}\neq {\cal O}_{M,x}$ for every $m >> 0$, 
by the classical theorem of Bombieri (\cite{b}). 
By Definition \ref{defAZD}, this implies that 
\[
\{ x \in X \mid n(\Theta_{h},x) > 0\} 
\subseteq \mbox{Supp}\,\mbox{Bs}\mid\!m_{0}L\!\mid
\]
holds. $\square$ 

\subsection{AZD for weakly pseudoeffective singular hermitian line bundles}
Similarly as Theorem \ref{AZD}, we obtain the following 
theorem. 
\begin{theorem}\label{AZD2}
Let $(L,h_{0})$ be a singular hermitian line bundle on a smooth projective variety $X$. 
Suppose that $(L,h_{0})$ is weakly pseudoeffetive.  Then  
\[
E(L,h_{0}):= \{ \varphi\in L^{1}_{loc}(X) \mid  \varphi \leqq 0, \,\,\Theta_{h_{0}} + \sqrt{-1}\partial\bar{\partial}\varphi \geqq 0 \}
\]
is nonempty and  if we define the function $\varphi_{P}\in L^{1}_{loc}(X)$ by
\[
\varphi_{P}(x)  := \sup \{\varphi (x)\mid \varphi \in E(L,h_{0})\} \,\,\,\, (x\in X).
\]
Then $h := e^{-\varphi_{P}}\cdot h_{0}$ is a singular hermitian metric on 
$L$ such that 
\begin{enumerate}
\item $\Theta_{h}\geqq 0$. 
\item $H^{0}(X,{\cal O}_{X}(mL)\otimes {\cal I}_{\infty}(h^{m})) 
\simeq H^{0}(X,{\cal O}_{X}(mL)\otimes {\cal I}_{\infty}(h_{0}^{m}))$ holds 
for every $m\geqq 0$. $\square$
\end{enumerate}
\end{theorem}
{\bf Proof of Theorem  \ref{AZD2}.}
Since $(L,h_{0})$ is pseudoeffective, there exists an ample line bundle 
$A$ on $X$ such that 
\[
H^{0}(X,{\cal O}_{X}(mL + A)\otimes {\cal I}(h_{0}^{m}))\neq 0
\]
holds for every $m\geqq 0$. Let $h_{A}$ be a $C^{\infty}$ hermitian metric 
on $A$ with strictly positive curvature and $dV$ be a $C^{\infty}$ volume form on $X$.  
Let 
\[
K_{m}:= K(mL +A,h_{0}^{m}\cdot h_{A},dV)
\]
be the (diagonal part) of the 
Bergman kernel of $mL +A$ with respect to the inner product 
\[
(\sigma ,\sigma^{\prime}) 
:= \int_{X}\sigma\cdot \bar{\sigma}^{\prime}\cdot h_{0}^{m}\cdot h_{A}\cdot dV
\]
on $H^{0}(X,{\cal O}_{X}(mL + A)\otimes {\cal I}(h_{0}^{m}))$. 

Let $x\in X$ be an arbitrary point and let 
$(U,z_{1},\cdot ,z_{n})$ be a coordinate neighbourhood centered at $x$ 
such that $U$ is biholomorphic to the open unit ball $B(O,1)$ in $\mathbb{C}^{n}$ centered at the origin via the coordinate.  
Taking $U$ to be sufficiently small, we may and do assume that 
there exist holomorphic frames $\mbox{\bf e}_{A},\mathbb{e}_{L}$ of $A$ and $L$ on $U$ respectively. 
Then with respect to these frames, we may express $h_{A},h_{L}$ as 
\[
h_{A} = e^{-\varphi_{A}}, h_{L} = e^{-\varphi_{L}}
\]
respectively in terms of plurisubharmonic functions $\varphi_{A},\varphi_{L}$
on $U$. 
By the extremal property of Bergman kernels, we see that 
$K_{m}(x) (x \in X)$ is expressed as :  
\[
K_{m}(x) = \sup \{ \mid\sigma (x)\mid^{2}\mid 
\sigma \in \Gamma (X,{\cal O}_{X}(A + mL)), \int_{X}\mid\sigma\mid^{2}\cdot h_{A}\cdot h_{L}^{m}\cdot dV = 1\}.
\]
Let $\sigma_{0} \in \Gamma (X,{\cal O}_{X}(A + mL))$ with 
\[
\int_{X}\mid\sigma_{0}\mid^{2}\cdot h_{A}\cdot h_{L}^{m}\cdot dV = 1
\]
and $\mid\sigma_{0}(x)\mid^{2} = K_{m}(x)$.
Let us write $\sigma_{0} = f\cdot\mbox{\bf e}_{A}\cdot\mbox{\bf e}_{L}^{m}$
on $U$ by using a holomorphic function $f$ on $U$.
By the submeanvalue property of plurisubharmonic functions, we have that 
\begin{eqnarray*}
\mid f(O)\mid^{2} &\leqq &  \frac{1}{\mbox{vol}(B(O,\varepsilon ))}\int_{B(O,\varepsilon)}\mid f\mid^{2}d\mu \\
&\leqq & (\sup_{B(O,\varepsilon )}e^{\varphi_{A}}\cdot e^{m\varphi_{L}})
\cdot (\frac{1}{\mbox{vol}(B(O,\varepsilon ))}\int_{B(,\varepsilon )}\mid f\mid^{2}e^{-\varphi_{A}}\cdot e^{-m\varphi_{L}}dV)\cdot 
(\sup_{B(O,\varepsilon )}\frac{d\mu}{dV})  \\
&\leqq & \frac{1}{\mbox{vol}(B(O,\varepsilon ))}\cdot (\sup_{B(O,\varepsilon )}e^{\varphi_{A}}\cdot e^{m\varphi_{L}})
\cdot (\sup_{B(O,\varepsilon )}\frac{d\mu}{dV})
\end{eqnarray*}
hold, where $d\mu$ is the standard Lebesgue measure on $\mathbb{C}^{n}$. 
Hence there exists a positive constant $C_{\varepsilon}$ independent of $m$
\[
K_{m}(x) \leqq C_{\varepsilon}\cdot \sup_{w\in B(O,\varepsilon )}(h_{A}^{-1}\cdot h_{L}^{-m})(w)\cdot dV
\] 
holds.  
Hence taking the $m$-th roots of the both sides, letting $m$ tend to infinity weand letting $\varepsilon$ tend to $0$,  we see that 
\[
K_{\infty}= \limsup_{m\rightarrow \infty}K(mL +A,h_{0}^{m}\cdot h_{A},dV)^{\frac{1}{m}}
\]
exists and satisfies the inequality 
\[
K_{\infty} \leqq h_{0}^{-1}  
\]
holds almost everywhere on $X$.  Hence if we set 
\[
h_{\infty} := \mbox{the lower envelope of}\,\, K_{\infty}^{-1} 
\]
is an element of $E(L,h_{0})$.   Hence $E(L,h_{0})$ is nonempty. 
The rest of the proof is the same as the one of Theorem \ref{AZD}. 
$\square$ 

\begin{definition}\label{singAZD}
Let $(L,h_{0})$ be a singular hermitian line bundle on a complex manifold $X$.
A singular hermitian metric $h$ on $L$ is said to be an analytic Zariski decomposition (AZD) of $(L,h_{0})$, if the followings hold : 
\begin{enumerate}
\item $\Theta_{h}\geqq 0$. 
\item $H^{0}(X,{\cal O}_{X}(mL)\otimes \bar{\cal I}_{\infty}(h^{m})) 
\simeq H^{0}(X,{\cal O}_{X}(mL)\otimes {\cal I}_{\infty}(h_{0}^{m}))$ holds 
for every $m\geqq 0$, where 
\[
\bar{\cal I}_{\infty}(h^{m}) := \cap_{p\geqq 1}\,\,{\cal I}_{p}(h^{m}). 
\]
 $\square$
\end{enumerate} 
\end{definition}
\begin{remark}
This definition is slightly different from that in \cite{tu7}. 
See Remark \ref{linfty}, for the reason why we use $L^{\infty}$
multiplier ideal sheaves  
instead of the usual multiplier ideal sheaves. $\square$  
\end{remark}
\begin{remark}\label{r2.4}
In Theorem \ref{AZD2}, $E(L,h_{0})$ is nonempty, if there exists a positive integer 
$m_{0}$ and $\sigma  \in  H^{0}(X,{\cal O}_{X}(m_{0}L)\otimes {\cal I}_{\infty}(h_{0}^{m_{0}}))$ such that $h_{0}^{m_{0}}(\sigma ,\sigma ) \leqq 1$. 
In this case 
\[
\varphi := \frac{1}{m_{0}}\log h_{0}^{m_{0}}(\sigma ,\sigma )
\]
belongs to $E(L,h_{0})$. $\square$ 
\end{remark}

About the normality of an AZD of a singular hermitian line bundle, 
we have the following proposition. 

\begin{proposition}\label{normality2}
Let $(L,h_{L})$ be a  pseudoeffective singluar hermitian line bundle on a smooth projective variety $X$ and  $F$ be a line bundle on $X$. 
Let $h_{F}$ be a $C^{\infty}$ hermitian metric on $F$.
 
Assume that $(L + F,h_{L}\cdot h_{F})$ is weakly pseudoeffective and let  $h$ be an AZD of $(L,h_{L})$. 

If $h_{L}$ is normal and 
\[
H^{0}(X,{\cal O}_{X}(m_{0}(L+F))\otimes {\cal I}_{\infty}(h_{L}^{m}))\neq 0
\]
, then $h$ is also normal. $\square$  
\end{proposition}

\subsection{$L^{2}$-extension theorem}
The $L^{2}$-extension theorem is our crucial tool to investigate 
multi adjoint bundles in this article.  
\begin{theorem}(\cite[p.200, Theorem]{o-t})\label{o-t}
Let $X$ be a Stein manifold of dimension $n$, $\psi$ a plurisubharmonic 
function on $X$ and $s$ a holomorphic function on $X$ such that $ds\neq 0$ 
on every branch of $s^{-1}(0)$.
We put $Y:= s^{-1}(0)$ and 
$Y_{0} := \{ x\in Y; ds(x)\neq 0\}$.
Let $g$ be a holomorphic $(n-1)$-form on $Y_{0}$ 
with 
\[
c_{n-1}\int_{Y_{0}}e^{-\psi}g\wedge\bar{g} < \infty ,
\]
where $c_{k}= (-1)^{\frac{k(k-1)}{2}}(\sqrt{-1})^{k}$. 
Then there exists a holomorphic $n$-form $G$ on 
$X$ such that 
\[
G(x) = g(x)\wedge ds(x) 
\]
on $Y_{0}$ and 
\[
c_{n}\int_{X}e^{-\psi}(1+\mid s\mid^{2})^{-2}G\wedge\bar{G} 
\leqq 1620\pi c_{n-1}\int_{Y_{0}}e^{-\psi}g\wedge\bar{g}. 
\] $\square$
\end{theorem}

\noindent For the extension from an arbitrary dimensional submanifold, 
T. Ohsawa extended Theorem \ref{o-t} in the following way. 

Let $X$ be a complex manifold of dimension $n$ and let $S$ be a closed complex submanifold of $X$. 
Then we consider a class of continuous function $\Psi : X\longrightarrow [-\infty , 0)$  such that  
\begin{enumerate}
\item $\Psi^{-1}(-\infty ) \supset S$,
\item if $S$ is $k$-dimensional around a point $x$, there exists a local 
coordinate $(z_{1},\ldots ,z_{n})$ on a neighbourhood of $x$ such that 
$z_{k+1} = \cdots = z_{n} = 0$ on $S\cap U$ and 
\[
\sup_{U\backslash S}\mid \Psi (z)-(n-k)\log\sum_{j=k+1}^{n}\mid z_{j}\mid^{2}\mid < \infty .
\]
\end{enumerate} 
The set of such functions $\Psi$ will be denoted by $\sharp (S)$. 

For each $\Psi \in \sharp (S)$, one can associate a positive measure 
$dV_{X}[\Psi ]$ on $S$ as the minimum element of the 
partially ordered set of positive measures $d\mu$ 
satisfying 
\[
\int_{S_{k}}f\, d\mu \geqq 
\limsup_{t\rightarrow\infty}\frac{2(n-k)}{v_{2n-2k-1}}
\int_{X}f\cdot e^{-\Psi}\cdot \chi_{R(\Psi ,t)}dV_{X}
\]
for any nonnegative continuous function $f$ with 
$\mbox{Supp}\, f\subset\subset X$.
Here $S_{k}$ denotes the $k$-dimensional component of $S$,
$v_{m}$ denotes the volume of the unit sphere 
in $\mathbb{R}^{m+1}$ and 
$\chi_{R(\Psi ,t)}$ denotes the characteristic funciton of the set 
\[
R(\Psi ,t) = \{ x\in M\mid -t-1 < \Psi (x) < -t\} .
\]

Let $X$ be a complex manifold and let $(E,h_{E})$ be a holomorphic hermitian vector 
bundle over $X$. 
Given a positive measure $d\mu_{X}$ on $X$,
we shall denote $A^{2}(X,E,h_{E},d\mu_{X})$ the space of 
$L^{2}$ holomorphic sections of $E$ over $X$ with respect to $h_{E}$ and 
$d\mu_{X}$. 
Let $S$ be a closed  complex submanifold of $X$ and let $d\mu_{S}$ 
be a positive measure on $S$. 
The measured submanifold $(S,d\mu_{S})$ is said to be a set of 
interpolation for $(E,h_{E},d\mu_{X})$, or for the 
sapce $A^{2}(X,E,h_{E},d\mu_{X})$, if there exists a bounded linear operator
\[
I : A^{2}(S,E\mid_{S},h_{E},d\mu_{S})\longrightarrow A^{2}(X,E,h_{E},d\mu_{X})
\]
such that $I(f)\mid_{S} = f$ for any $f \in  A^{2}(S,E\mid_{S},h_{E},d\mu_{S})$. 
$I$ is called an interpolation operator.
The following theorem is crucial.

\begin{theorem}(\cite[Theorem 4]{o})\label{extension}
Let $X$ be a complex manifold with a continuous volume form $dV_{X}$,
let $E$ be a holomorphic vector bundle over $X$ with $C^{\infty}$ fiber 
metric $h_{E}$, let $S$ be a closed complex submanifold of $X$,
let $\Psi\in \sharp (S)$ and let $K_{X}$ be the canonical bundle of $X$.
Then $(S,dV_{X}(\Psi ))$ is a set of interpolation 
for $(E\otimes K_{X},h_{E}\otimes (dV_{X})^{-1},dV_{X})$, if 
the followings are satisfied.
\begin{enumerate}
\item There exists a closed set $F\subset X$ such that 
\begin{enumerate}
\item $F$ is locally negligble with respect to $L^{2}$-holomorphic functions, i.e., 
for any local coordinate neighbourhood $U\subset M$ and for any $L^{2}$-holomorphic function $f$ on $U\backslash X$, there exists a holomorphic function 
$\tilde{f}$ on $U$ such that $\tilde{f}\mid U\backslash F = f$.
\item $M\backslash F$ is a Stein manifold which intersects with every component of $S$. 
\end{enumerate}
\item $\Theta_{h_{E}}\geqq 0$ in the sense of Nakano,
\item $\Psi \in \sharp (S)\cap C^{\infty}(X\backslash S)$,
\item $e^{-(1+\epsilon )\Psi}\cdot h_{E}$ has semipositive 
curvature in the sense of Nakano for every $\epsilon \in [0,\delta]$ 
for some $\delta > 0$.
\end{enumerate}
Under these conditions, there exists a constant $C$ and an interpolation operator 
from $A^{2}(S,E\otimes K_{X}\mid_{S},h\otimes (dV_{X})^{-1}\mid_{S},dV_{X}[\Psi ])$
to \\ $A^{2}(X,E\otimes K_{X},h\otimes (dV_{X})^{-1}.dV_{X})$ whose 
norm does not exceed $C\delta^{-3/2}$.
If $\Psi$ is plurisubharmonic, the interpolation operator can be chosen 
so that its norm is less than $2^{4}\pi^{1/2}$. $\square$
\end{theorem}

The above theorem can be generalized to the case that 
$(E,h_{E})$ is a singular hermitian line bundle with semipositive
curvature current  (we call such a singular hermitian line 
bundle $(E,h_{E})$ a {\bf pseudoeffective singular hermitian line bundle}) as was remarked in \cite{o}. 

\begin{lemma}\label{extension2}
Let $X,S,\Psi ,dV_{X}, dV_{X}[\Psi], (E,h_{E})$ be as in Theorem \ref{extension} 
Let $(L,h_{L})$ be a pseudoeffective singular hermitian line 
bundle on $X$. 
Then $S$ is a set of interpolation for 
$(K_{X}\otimes E\otimes L,dV_{X}^{-1}\otimes h_{E}\otimes h_{L})$.  $\square$
\end{lemma}

Later we shall use the more general residue volume form $dV[\Psi ]$ 
as is introduced in Section 1.2. 
Since the proof of Theorem \ref{extension} in \cite{o} works without any 
change, we shall also apply Theorem \ref{extension} also for this generalized 
residue volume form $dV[\Psi ]$, too.

\section{Numerical Kodaira dimension of pseudoeffective singular hermitian line bundles}

Let $X$ be a smooth projective variety and let $L$ be a line bundle on $X$.
\begin{definition}\label{L-dim}(\cite{nak})
The {\bf $L$-dimension}  $\mbox{\em Kod} (L)$ is defined by 
\[
\mbox{\em Kod} (L): = \limsup_{m\rightarrow\infty}\frac{\log \dim H^{0}(X,{\cal O}_{X}(mL))}
{\log m}. 
\]
The {\bf numerical dimension} $\nu (L)$ of $L$  is defined by 
\[
\nu (L):= \sup_{A}\limsup_{m\rightarrow\infty}\frac{\log \dim H^{0}(X,{\cal O}_{X}(A +mL))}{\log m},
\] 
where $A$ runs all the ample line bundles on $X$. 
We often denote $\mbox{\em Kod} (K_{X})$ by $\mbox{\em Kod} (X)$ and call it 
the {\bf Kodaira dimension} of $X$.  And we often denote $\nu (K_{X})$ by 
$\nu (X)$ and call it the {\bf numerical Kodaira dimension} of $X$.  
$\square$ 
\end{definition} 
It is trivial to see that $\mbox{Kod} (L)\leqq \nu (L)$ holds. 
And $\mbox{Kod} (L),\nu (L)$ are either $-\infty$ or integers 
between $0$ and $\dim X$. 

The purpose of this section is to define a similar numerical Kodaira dimension for 
a psuedoeffective singular hermitian line bundle on a smooth projective 
variety and study the relation between a numerical property of the singular 
hermitian line bundle and the asymptotic property of the powers of the singular hermitian line bundle twisted by a sufficiently positive line bundle.  

\subsection{Intersection theory for pseudoeffective singular hermitian 
line bundles}

To intoroduce the notion of  the numerical Kodaira dimension of a pseudoeffective singular hermitian 
line bundle on a smooth projective variety,  first we define the intersection 
number.  
\begin{definition}(\cite{tu})\label{rmwp}
Let $(L,h_{L})$ be a weakly pseudoeffective singular hermitian line bundle on 
a smooth projective $n$-fold $X$. 
 The intersection number 
$(L,h_{L})^{n}\cdot X$ defined by 
\[
(L,h_{L})^{n}\cdot X := n!\cdot \limsup_{m\rightarrow\infty}m^{-n}
\dim H^{0}(X,{\cal O}_{X}(mL)\otimes {\cal I}(h^{m})).
\]
$(L,h_{L})$ is said to be big, if $(L,h_{L})^{n}\cdot X$ is positive. 
For a $r$-dimensional subvariety $V$ in $X$ such that 
$h_{L}\mid_{V}$ is not identically $+\infty$, we define 
\[
(L,h_{L})^{r}\cdot V:=  r!\cdot \limsup_{m\rightarrow\infty}m^{-r}
\dim H^{0}(\tilde{V},{\cal O}_{\tilde{V}}(m\mu^{*}L)\otimes {\cal I}(\mu^{*}(h_{L}\!\!\mid_{V})^{m})), 
\] 
where $\mu : \tilde{V} \longrightarrow V$ is a resolution of 
singularities.  
$(L,h_{L})$ is said to be big on $V$, if $(L,h_{L})^{r}\cdot V > 0(r = \dim V)$ 
holds. 
 $\square$ 
\end{definition}
The well definedness of the intersection number is verified as follows. 
\begin{proposition}\label{well defined}
The definition of $(L,h_{L})^{r}\cdot V$ is independent of the choice 
of the resolution $\pi : \tilde{V} \longrightarrow V$. $\square$ 
\end{proposition}
{\bf Proof.}
Let $\mu : \tilde{V} \longrightarrow V$ be an resolution and let 
$\mu^{\prime} : \tilde{V}^{\prime}$ be another resolution 
factors through $\mu$, i.e., there exists a morphism 
$\phi : \tilde{V}^{\prime} \longrightarrow \tilde{V}$ such that 
$\mu^{\prime} = \mu\circ \phi$. 
Then 
\[
\phi_{*}({\cal O}_{\tilde{V}^{\prime}}(K_{\tilde{V}^{\prime}})\otimes {\cal I}((\mu^{\prime})^{*}(h\mid_{V})^{m}))) = {\cal O}_{V}(K_{V})\otimes {\cal I}(\mu^{*}(h\mid_{V})^{m}))
\]
holds.  Hence 
\begin{equation}\label{blow up}
H^{0}(\tilde{V},{\cal O}_{\tilde{V}}(K_{\tilde{V}}+ m\mu^{*}L)\otimes {\cal I}(\mu^{*}(h\mid_{V})^{m})) = H^{0}(\tilde{V^{\prime}},{\cal O}_{\tilde{V^{\prime}}}(K_{\tilde{V}^{\prime}}+ m(\mu^{\prime})^{*}L)\otimes {\cal I}((\mu^{\prime})^{*}(h\mid_{V})^{m})) 
\end{equation}
holds for every $m\geqq 1$. 
On the other hand 
\[
\limsup_{m\rightarrow\infty}m^{-r}\dim 
H^{0}(\tilde{V},{\cal O}_{\tilde{V}}(K_{\tilde{V}}+ m\mu^{*}L)\otimes {\cal I}(\mu^{*}(h\mid_{V})^{m})) 
=  \vspace{-10mm}
\]
\[
\hspace{30mm}
\limsup_{m\rightarrow\infty}m^{-r}\dim 
H^{0}(\tilde{V},{\cal O}_{\tilde{V}}(m\mu^{*}L)\otimes {\cal I}(\mu^{*}(h\mid_{V})^{m})) 
\]
holds, since 
\[
\limsup_{m\rightarrow\infty}m^{-r}\dim 
H^{0}(\tilde{V},{\cal O}_{\tilde{V}}(m\mu^{*}L + A - B)\otimes {\cal I}(\mu^{*}(h\mid_{V})^{m}))
\]
\[
\hspace{30mm} 
= 
\limsup_{m\rightarrow\infty}m^{-r}\dim 
H^{0}(\tilde{V},{\cal O}_{\tilde{V}}(m\mu^{*}L)\otimes {\cal I}(\mu^{*}(h\mid_{V})^{m})) 
\]
holds for any very ample divisors $A$ and $B$ on $\tilde{V}$.
In fact this can be verified by using the exact sequence
\[
0 \rightarrow 
H^{0}(\tilde{V},{\cal O}_{\tilde{V}}(m\mu^{*}L + A)\otimes {\cal I}(\mu^{*}(h\mid_{V})^{m}))
\rightarrow H^{0}(\tilde{V},{\cal O}_{\tilde{V}}(m\mu^{*}L + A - B)\otimes {\cal I}(\mu^{*}(h\mid_{V})^{m}))
\]
\[
\hspace{50mm} 
\rightarrow 
H^{0}(B,{\cal O}_{\tilde{V}}(m\mu^{*}L + A))
\] 
, etc. 
And similar equality holds on $\tilde{V}^{\prime}$. 
Hence using the equality (\ref{blow up}), we conclude that 
\[
\limsup_{m\rightarrow\infty}m^{-r}\dim 
H^{0}(\tilde{V},{\cal O}_{\tilde{V}}(m\mu^{*}L)\otimes {\cal I}(\mu^{*}(h\mid_{V})^{m})) 
\]
\[
\hspace{30mm} 
= \limsup_{m\rightarrow\infty}m^{-r}\dim 
H^{0}(\tilde{V}^{\prime},{\cal O}_{\tilde{V}^{\prime}}(m\mu^{*}L)\otimes {\cal I}((\mu^{\prime})^{*}(h\mid_{V})^{m})) 
\]
holds.   This completes the proof of Proposition \ref{well defined}. $\square$

\subsection{The numerical Kodaira dimension and the asymptotic Kodaira dimension of  a pseudoeffective singular hermitian line bundle}\label{subsec}
 
We shall define the numerical Kodaira dimension  and the asymptotic 
Kodaira dimension for pseudoeffective 
singular hermitian line bundles on smooth projective varieties.
These invariants play the essential roles in this article.  

\begin{definition}\label{numerical}
$(L,h_{L})$ be a pseudoeffective singular hermitian line bundle on a projective  manifold  $X$.  
We set 
\[
\nu_{num}(L,h_{L}) := \sup \{\dim V \mid \mbox{\em $V$ is a subvariety of $X$ such that 
$h_{L}\!\mid_{V}$ is well defined} 
\]
\[
\hspace{30mm}\mbox{\em and $(L,h_{L})^{\dim V}\!\!\cdot V> 0$}\}. 
\]  
We call $\nu_{num}(L,h_{L})$ the {\bf numerical Kodaira dimension} of $(L,h_{L})$.
  $\square$
\end{definition}
The following invariant is more analytic in nature. 
\begin{definition}\label{numerical Kodaira}
The {\bf asymptotic Kodaira dimension} $\nu (L,h_{L})$ of $(L,h_{L})$ is defined by 
\[
\nu_{asym} (L,h_{L}) = \sup_{A}\limsup_{m\rightarrow\infty}\frac{\log h^{0}(X,{\cal O}_{X}(A + mL)\otimes {\cal I}(h_{L}^{m}))}{\log m}, 
\]
where $A$ runs all the ample line bundles on $X$. $\square$
\end{definition}
The following example shows that $\nu_{asym}$ is not necessarily an 
integer.  
\begin{example}\label{nonintegral}
Let $T$ be a closed positive $(1,1)$ current on $\mathbb{P}^{1}$
\[
T =  \sum_{i =1}^{\infty}\sum_{j=1}^{3^{n+1}}\frac{1}{4^{n}}P_{ij}
\]
where $\{P_{ij}\}$ are distinct points on $\mathbb{P}^{1}$. 
Then there exists a singular hermitian metric $h$ on ${\cal O}(1)$
such that $\Theta_{h} = 2\pi T$. 
Then we see that 
\[
\nu_{asym}({\cal O}(1),h) = \frac{\log 3}{\log 4} 
\]
and 
\[
\nu_{num}({\cal O}(1),h) = 0.
\]
This implies that $\nu_{num} \neq \nu_{asym}$ in general. 
$\square$
\end{example}

\subsection{Seshadri constant for a pseudoeffective singular hermitian line bundle}

In this subsection we shall give a criterion of the bigness of a 
psueodeffective singular hermitian line bundles on smooth projective varieties. 
\begin{definition}
Let $X$ be a smooth projective variety and let $H$ be an ample divisor on $X$.
Let $(L,h_{L})$ be a 
psuedoeffective singular hermitian line bundle on $X$.  
Let $x$ be a point on $X$. 
We set  
\[
\epsilon ((L,h_{L}),H,x) = \inf_{C}\frac{(L,h_{L})\cdot C}{H\cdot C}, 
\]
where $C$ runs all the irreducible curves in $X$ passing through $x$. 
We call $\epsilon ((L,h_{L}),H,x)$  the Seshadri constant of $(L,h_{L})$ at 
$x$ with respect to $H$. $\square$
\end{definition}

\begin{theorem}\label{nonbig}
Let $X$ be a smooth projective variety  and let $H$ be an ample divisor on $X$.
Let  $(L,h_{L})$ be a 
pseudoeffective singular hermitian line bundle on $X$.

Then there exists at most  countable union of proper subvarieties $F$ such that 
if  there exists a point $x_{0}\in X - F$ such that 
\[
\epsilon ((L,h_{L}),x_{0}) > 0
\]
holds, then $(L,h_{L})$ is big.  $\square$
\end{theorem}
{\bf Proof of Theorem \ref{nonbig}.}
We say that $C$ is a strongly movable curve, if 
\[
C = \mu_{*}(\tilde{A}_{1}\cap \cdots \cap \tilde{A}_{n-1})
\]
for some very ample divisors $\tilde{A}_{j}$ on $\tilde{X}$, where 
$\mu :\tilde{X}\longrightarrow X$ is a modification. 
The strongly movavble cone $SME(X)$ of $X$ is the cone of curves   
generated by all the strongly movable curves on $X$. 
Let ${\cal S}$ be the family of strongly movable curves on $X$.  
We set 
\[
U  := \{ x\in X\mid  \mbox{For every irreducible component of ${\cal C}$, there exists an irreducible member $C$  }
\]
\[
\hspace{30mm}
\mbox{
belonging to the component and 
passing through $x$.}\}.
\]
Then we see that there exists at most a coutable union $F$ of proper subvarieties of $X$ such that $U = X - F$. 

Suppose that there exists a point $x_{0}\in U$ such that 
$\epsilon_{0}:= \epsilon ((L,h_{L}),H,x) > 0$. 
Let $dV$ be a $C^{\infty}$ volume form on $X$. 
Let  $m$ be a positive integer such that $m > 2/\epsilon_{0}$
and  let 
\[
H^{0}(X,{\cal O}_{X}(m_{0}L + H)\otimes {\cal I}(h^{m_{0}}_{L})) 
\]
us consider the inner product 
\[
(\sigma,\sigma^{\prime}):= \int_{X}\sigma\cdot \bar{\sigma}^{\prime}\cdot h_{L}^{m}\cdot h_{H}\cdot dV. 
\]
Let $K_{m}:= K(X,mL + H,h_{L}^{m}\cdot h_{H},dV)$ be the diagonal part of 
the Bergman kernel with respect to the above inner product.
We set 
\[
h_{m}:= (K_{m})^{-\frac{1}{m}}. 
\]
Then $h_{m}$ is a singular hermitian metric on $L + \frac{1}{m}H$ 
with algebraic singularities.  

We note that for every strongly movable 
curve $C$ passing through $x_{0}$ 
\[
(L +  \frac{1}{m}H,h_{m}\mid C)\cdot C \geq  \frac{\varepsilon_{0}}{2}H\cdot C 
\]
holds.
Since the pseudoeffective cone of $X$ is the dual of $SME(X)$ as in \cite{b-d-p-p}, by the definition of $U$,  we see that 
\[
(L+\frac{1}{m}H,h_{m}) - \frac{\epsilon_{0}}{2}H
\]
is pseudoeffective.  
Letting $m$ tend to infinity, by Lemma \ref{reciprocity}, we see that 
\[
(L- \frac{\epsilon_{0}}{2}H,h_{L}\cdot h_{H}^{-\frac{\epsilon_{0}}{2}})
\]
is pseudoeffective (cf. Definition \ref{pe}).  
Hence $(L,h_{L})$ should be big.  
This completes the proof of Theorem \ref{nonbig}.  $\square$ 

\subsection{Non big pseudoeffective singular hermitian line bundles}

In this subsection, we shall prove the following vanishing theorem. 
\begin{proposition}\label{vanish}
Let $(L,h_{L})$ be a pseudoeffective singular hermitian line bundle
on a smooth projective variety $X$.  Suppose that $(L,h_{L})$ is not big and  
one of the followings holds. 
\begin{enumerate}
\item $L$ is not big.
\item  $L$ is normal (cf. Definition \ref{normal}).
\end{enumerate}
Then there exists a very ample divisor $H$ on 
$X$  such that 
\[
H^{0}(X,{\cal O}_{X}(mL - H)\otimes {\cal I}(h_{L}^{m})) \neq 0 
\] 
holds for every $m \geqq 1$. $\square$  
\end{proposition}
{\bf Proof of Proposition \ref{vanish}.}
Suppose that there exists a very ample divisor $H$ on $X$ such that 
\[
H^{0}(X,{\cal O}_{X}(mL - H)\otimes {\cal I}(h_{L}^{m})) \neq 0 
\] 
holds for some  $m \geqq 1$. 
Then since $H$ is very ample, we see that 
\[
\mid H^{0}(X,{\cal O}_{X}(mL)\otimes {\cal I}(h_{L}^{m}))\mid 
\]
gives a birational rational map from $X$ into a projective space. 
Hence if $L$ is not big, then for every very ample divor $H$,
\[
H^{0}(X,{\cal O}_{X}(mL - H)\otimes {\cal I}(h_{L}^{m})) = 0 
\] 
holds for every $m \geqq 1$.  \vspace{3mm}\\

\noindent Next we shall assume that $(L,h_{L})$ is normal. 
If we take $H$ very general, we see that 
\[
{\cal I}(h_{L}^{m}\mid_{H}) = {\cal I}(h_{L}^{m})\mid_{H}
\]
holds for every $m\geqq 1$.
This is possible by the following lemma. 
\begin{lemma}\label{generic}
There exists a smooth member $H^{\prime}\in \mid\! H\!\mid$, such that 
\[
{\cal I}(h_{L}^{m})\otimes{\cal O}_{H^{\prime}}
= {\cal I}(h_{L}^{m}\mid_{H^{\prime}})
\]
holds for every $m\geqq 1$. $\square$   
\end{lemma}
{\bf Proof of Lemma \ref{generic}}.
Let $A$ be a sufficiently ample line bundle such that 
${\cal O}_{X}(A+mL))\otimes {\cal I}(h_{L}^{m})$  is globally generated
for all $m\geqq 1$.  
Let $\{ \sigma^{(m)}_{j}\}_{j=1}^{N_{m}}$ be a (complete) basis 
of $H^{0}(X,{\cal O}_{X}(A+mL))\otimes {\cal I}(h_{L}^{m}))$. 
We consider the subset
\[
U:= \{ F\in \mid\! H\!\mid ; \mbox{$F$ is smooth}, 
\int_{F}\mid\sigma_{j}^{(m)}\mid^{2}\cdot h_{L}^{m}\cdot h_{A}\cdot \, dV_{F}
< + \infty \,\,\,
\]
\vspace{-8mm}
\[
\hspace{50mm}\mbox{for every $m$ and $1\leqq j\leqq N_{m}$.}\}
\]
of $\mid\! H\!\mid$, where $dV_{F}$ denotes the volume form  on $F$ induced by the 
K\"{a}hler form $\omega$. 
We claim that such $U$ is the complement of at most a countable union of 
proper subvarieties of $\mid\! H\!\mid$.
Let us fix a positive integer $m$.  
Since ${\cal I}(h_{L}^{m})$ is a coherent sheaf (\cite{n}), we see that
\[
U_{m}:= \{ F\in \mid\! H\!\mid ; \mbox{$F$ is smooth}, 
\int_{F}\mid\sigma_{j}^{(m)}\mid^{2}\cdot h_{L}^{m}\cdot h_{A}\,\, dV_{F}
< + \infty \,\,\,
\] 
\vspace{-8mm}
\[
\hspace{50mm}\mbox{for $1\leqq j\leqq N_{m}$.}\}
\] 
is a Zariski open subset of $\mid\! H\!\mid$.
In fact this can be verified as follows. 
Let $\Lambda$ be a pencil contained in $\mid\! H\!\mid$ which contains a 
smooth member. 
Then by Fubini's theorem, we see that a general member of 
$\Lambda$ is contained in $U_{m}$, unless for every general member $F$ 
of $\Lambda$ the set 
\[
\{x \in F\mid h_{A}\cdot h_{L}^{m}\mid\sigma^{(m)}_{j}\mid^{2}
\not{\in} L^{1}_{loc}(F,x) \,\,\mbox{for some $1\leqq j\leqq N_{m}$}\}
\]
is contained in the base locus of $\Lambda$. 
Hence we see that $U_{m}$ is Zariski dense in $\mid\! H\!\mid$. 
Then since $U = \cap_{m=1}^{\infty}U_{m}$, we complete the proof of 
Lemma \ref{generic}.   
$\square$.  \vspace{5mm} \\
Let us continue the proof of Proposition \ref{vanish}. 
Since $(L,h_{L})$ is not big, by Theorem \ref{nonbig} we see that 
for a very general $x\in X$, $\epsilon ((L,h_{L}),H,x) = 0$ holds.
\begin{lemma}\label{very general}
There exists at most a countable union  of subvarities $F$ 
 such that for every $x\in  X - F$ and $\delta > 0$, there exists an irreducible 
curve $C$ on $X$  passing through $x$ such that
\begin{enumerate}
\item $C$ is smooth at $x$, 
\item $(L,h_{L})\cdot C \leqq \delta\cdot (H\cdot C)$. 
\end{enumerate}  
$\square$ 
\end{lemma}
By the $L^{2}$-extension theorem  (Theorems \ref{o-t} and \ref{extension}) the intersection number $(L,h_{L})\cdot C_{t}$ is
lower semicontinuous with respect to the countable Zariski topology 
on $\Delta^{n-1}$.  
Then since the Hilbert scheme (of $X$) has only countably many irreducible components, by Lemma \ref{very general},  we have the following lemma. 
\begin{lemma} \label{general2}
For every $x\in X - F$ and $\delta > 0$ there exist a 
$(n-1)$-dimensional family of irreducible curves $\{ C_{t}\}_{t\in\Delta^{n-1}}$parametrized 
$\Delta^{n-1}$ and  a coordinate neighbourhood $(U, z_{1},\cdots ,z_{n})$ of $x$  such that 
\begin{enumerate}
\item $(L,h_{L})\cdot C_{t} \leqq \delta (H\cdot C_{t})$ holds for every 
$t\in \Delta^{n-1}$,@ 
\item $z_{1}(x)= \cdots = z_{n}(x) = 0$,
\item $U$ is biholomorphic to $\Delta^{n}$ via the coordinate 
$(z_{1},\cdots ,z_{n})$. 
\item For every $t\in \Delta^{n-1}$, 
\[
C_{t} \cap U = \{ p\in U \mid (z_{2}(p),\cdots ,z_{n}(p)) = t\}
\]
holds for every $t\in \Delta^{n-1}$. 
\end{enumerate}
$\square$  
\end{lemma}
By the assumption $E$ is contained in a proper analytic subset of $X$, say $V$. 
Taking $H$ to be sufficiently ample, we may and do assume that 
$\frac{1}{2}H$ is Cartier and $V$ is 
contained in a member $H_{0}$ of $\mid \frac{1}{2}H\mid$. 
Let us take $\delta < 1/2m$ in  Lemma \ref{general2}. 
Then  since 
\[
\deg {\cal O}_{C_{t}}(mL -  H) \otimes {\cal I}(h_{L}^{m}\mid_{C_{t}})
< - \frac{1}{2m}\cdot m \cdot(H\cdot C) + \frac{1}{2}H\cdot C < 0    
\] 
hold (because $V$ is 
contained in a member $H_{0}\in \mid\frac{1}{2}H\mid$), we see that 
\[
H^{0}(C_{t},{\cal O}_{C_{t}}(mL -  H) \otimes {\cal I}(h_{L}^{m}\mid_{C_{t}})) = 0
\]
holds for every $t \in \Delta^{n-1}$. 
By Fubini's theorem and Lemma \ref{general2},
\[
H^{0}(X,{\cal O}_{X}(mL - H)\otimes
{\cal I}(h^{m}_{L})) = 0
\]
holds.  
We note that $H$ can be taken independent of $m$. 
This completes the proof of Proposition \ref{vanish}. $\square$ 
   
\subsection{Relation between $\nu_{num}$ and $\nu_{asym}$}
In Example \ref{nonintegral}, we have seen that  
$\nu_{num}$ and $\nu_{asym}$ are different in general. 
But in this subsection we shall prove that $\nu_{num} = \nu_{asym}$ holds 
under a mild condition. 

\begin{theorem}\label{nak}
Let $X$ be a smooth projective variety and let $(L,h_{L})$ be a 
pseudoeffecive singular hermitian line bundle on $X$.

Then 
\[
\nu_{num}(L,h_{L}) \leqq  \nu_{asym}(L,h_{L})
\]
holds. 
Moreover if  $h_{L}$ is normal (cf. Definition \ref{normal}),
for every ample line bundle $A$ on $X$,   
\[
\dim H^{0}(X,{\cal O}_{X}(A + mL)\otimes{\cal I}(h_{L}^{m}))
= O(m^{\nu_{num}(L,h_{L}))}).
\]
holds. 
In particular 
\[
\nu_{num}(L,h_{L}) =  \nu_{asym}(L,h_{L})
\]
holds. 
$\square$
\end{theorem}
\begin{remark}
If $h_{L}$ is an AZD of $L$ with $\mbox{Kod}(L)\geqq 0$, then 
the set 
\[
E := \{ x\in X\mid n(\Theta_{h_{L}},x) > 0\}
\]
is contained in the stable base locus of 
$L$.  Hence in this case $h_{L}$ is normal.   
$\square$ 
\end{remark}
\noindent {\bf Proof of Theorem \ref{nak}}.
We denote $\nu_{num}(L,h_{L})$ by $\nu$.  
Let $V$ be a $\nu$ dimensional subvariety such that $h_{L}\mid_{V}$ is well defined 
and $(L\!\!\mid_{V},h_{L}\!\!\mid{V})$ is big.
Let $f : Y \longrightarrow X$ be an embedded resolution of $V$.
Then replacing $(L,h_{L})$ by $f^{*}(L,h_{L})$ we may assume that 
$V$ is smooth from the beginning.  
Let $A$ be a sufficiently ample line bundle such that 
every element of 
\[
H^{0}(V,{\cal O}_{V}(A + mL)\otimes {\cal I}(h_{L}^{m}\mid_{V}))
\]
extends to an element $H^{0}(X,{\cal O}_{X}(A + mL)\otimes {\cal I}(h_{L}^{m}))$. 
Then we have that 
\begin{equation}\label{leq}
\nu_{num}(L,h_{L}) \leqq  \limsup_{m\rightarrow\infty}\frac{\log \dim H^{0}(X,{\cal O}_{X}(A + mL)
\otimes {\cal I}(h^{m}_{L}))}{\log m}
\end{equation}
holds.  Hence we have the inequality
\[
\nu_{num}(L,h_{L}) \leqq \nu_{asym}(L,h_{L}).
\]

Next suppose that  $(L,h_{L})$ be a normal  pseudoeffective  singular hermitian line bundle on 
a smooth projective variety $X$. 
We shall prove 
\[
\nu_{num} (L,h_{L}) \geqq \limsup_{m\rightarrow\infty}
\frac{\log \dim H^{0}(X,{\cal O}_{X}(A+ mL)\otimes {\cal I}(h^{m}_{L}))}
{\log m}
\]
holds by induction on $n= \dim X$.

If $n= 1$ and  $(L,h_{L})$ is not big, then $(L,h_{L})$ is numerically trivial.
Hence  $\Theta_{h_{L}}$ has no absolutely continuous part. 
Since $h_{L}$ is normal,   
\[
\Theta_{h_{L}} = \sum a_{i}P_{i}
\]
for some effective $\mathbb{R}$-divisor $\sum a_{i}P_{i}$ on $X$. 
Hence $\nu_{num}(L,h_{L}) = \nu_{asym}(L,h_{L}) = 0$ holds in this case. 

Suppose that 
\[
\nu_{num}(F,h_{F}) =  \limsup_{m\rightarrow\infty}
\frac{\log \dim H^{0}(Y,{\cal O}_{F}(A+ mF)\otimes {\cal I}(h^{m}_{F}))}
{\log m}
\]
holds for every normal pseudoeffective singular hermitian line bundle $(F,h_{F})$ and every sufficiently ample line bundle $A$ on a smooth projective variety $Y$ of dimension $ \leqq n-1$.

Let $X$ be a smooth projective variety of dimension $n$ and let $(L,h_{L})$ be 
a normal  pseudoeffective line bundle on $X$.  
If $(L,h_{L})$ is big, there is nothing to prove. 
We shall assume  that $(L,h_{L})$ is not big, i.e., 
\[
\limsup_{m\rightarrow\infty}m^{-n}\dim H^{0}(X,{\cal O}_{X}(mL)\otimes 
{\cal I}(h^{m})) = 0
\]
holds. 
Let $H$ be a sufficiently ample very ample smooth divisor on $X$.
Then by Proposition \ref{vanish}, 
\begin{equation}\label{va}
H^{0}(X,{\cal O}_{X}(mL -H)\otimes {\cal I}(h_{L}^{m})) = 0
\end{equation}
holds for every $m \geqq 0$.
Let $G$ be a smooth member of $\mid\!2H\!\mid$.   
If we take $G$ properly,by the Lemma \ref{generic} we may assume that 
\[
{\cal I}(h_{L}^{m})\otimes {\cal O}_{G} 
= {\cal I}(h_{L}^{m}\mid_{G})
\]
holds for every $m \geqq 1$.    
Let us consider the exact sequence
\[
0 \rightarrow H^{0}(X,{\cal O}_{X}(mL - H)\otimes {\cal I}(h_{L}^{m}))
\rightarrow H^{0}(X,{\cal O}_{X}(mL + H)\otimes {\cal I}(h_{L}^{m}))
\]
\[
\hspace{60mm} \rightarrow H^{0}(G,{\cal O}_{G}(mL + H)\otimes {\cal I}(h_{L}^{m})). 
 \]
 Then by (\ref{va}) we have that 
\[
\dim H^{0}(X,{\cal O}_{X}(mL +H)\otimes {\cal I}(h^{m}_{L}))
\leqq 
\dim H^{0}(G,{\cal O}_{G}(mL +H)\otimes {\cal I}(h^{m}_{L})).
\]
Since ${\cal I}(h_{L}^{m}\mid_{G}) = {\cal I}(h_{L}^{m})\mid_{G}$ holds 
for every $m \geqq 0$ by the choice of $G$, we see that  
\begin{equation}\label{asymp}
\limsup_{m\rightarrow\infty}
\frac{\log \dim H^{0}(X,{\cal O}_{X}(H + mL)\otimes {\cal I}(h^{m}_{L}))}{\log m}
\leqq  \limsup_{m\rightarrow\infty}
\frac{\log \dim H^{0}(G,{\cal O}_{G}(mL +H)\otimes {\cal I}(h^{m}_{L}\mid G))}{\log m}
\end{equation}
holds.  By the induction assumption 
\begin{equation}\label{induction}
\limsup_{m\rightarrow\infty}\frac{\log \dim H^{0}(G,{\cal O}_{G}(mL +H)\otimes {\cal I}(h^{m}_{L}\mid_{G}))}{\log m}
=   
\nu_{num}(L\mid_{G},h_{L}\mid_{G})
\end{equation}
holds.  By the assumption, if we take $H$ sufficiently ample, we see that 
\begin{equation}\label{equality}
\nu_{num}(L,h_{L}) = \nu_{num}(L\mid_{G},h_{L}\mid_{G})
\end{equation}
holds. 
Hence combining  (\ref{asymp}),(\ref{induction}),(\ref{equality}), we see that 
\begin{equation}\label{geq}
\nu_{num}(L,h_{L})\geqq \limsup_{m\rightarrow\infty}
\frac{\log \dim H^{0}(X,{\cal O}_{X}(H + mL)\otimes {\cal I}(h^{m}_{L}))}{\log m}
\end{equation}
holds. 

Combining (\ref{leq}) and (\ref{geq}), we see that 
\[
\nu_{num}(L,h_{L}) = \nu_{asym}(L,h_{L})
\]
holds.  

By using the construction as above inductively, we find a sequence 
of members $G = G_{1}, G_{2}, \cdots , G_{n-\nu} \in \mid\!2H\!\mid (n:= \dim X, \nu = \nu_{num}(L,h_{L}))$ such that 
\begin{enumerate}
\item $G_{i}$ intersects $G_{1}\cap \cdots \cap G_{i-1}$ transversally.
We set $X_{i}:= G_{1}\cap\cdots \cap G_{i}$ for $1\leqq i\leqq n-\nu$ 
and $X_{0} := X$.  
\item $h_{L}\!\mid_{X_{i}}$ is well defined for every $1\leqq i\leqq n-\nu$. 
\item ${\cal I}(h_{L}^{m})\mid_{X_{i}} = {\cal I}(h^{m}_{L}\mid_{X_{i}})$ 
for every $1\leqq i\leqq n-\nu$.
\item 
\[
\dim H^{0}(X_{i},{\cal O}_{X_{i}}(H +mL)\otimes {\cal I}(h^{m}_{L}))
\leqq \dim H^{0}(X_{i+1},{\cal O}_{X_{i+1}}(H+mL)\otimes {\cal I}(h^{m}_{L}))
\]
holds for every $m\geqq 1$ and $0\leqq i\leqq n-\nu-1$.
\item $(L,h_{L})\!\mid_{X_{n-\nu}}$ is big.  \vspace{3mm}
\end{enumerate}
Then we see that 
\[
\dim H^{0}(X,{\cal O}_{X}(H + mL)\otimes {\cal I}(h_{L}^{m}))
\leqq \dim H^{0}(X_{n-\nu},{\cal O}_{X_{\nu}}(H + mL)\otimes {\cal I}(h_{L}^{m}\mid_{X_{n-\nu}}))
\]
holds.  
Since 
\[
\dim H^{0}(X_{n-\nu},{\cal O}_{X_{\nu}}(H + mL)\otimes {\cal I}(h_{L}^{m}\mid_{X_{n-\nu}})) = O(m^{\nu}),  
\]
this completes the proof of 
Theorem \ref{nak}. $\square$

\section{Asymptotic expansion of Bergman kernels of pseudoeffective
line bundles}

The purpose of this section is to extend the asymptotic expansions of 
Bergman kernels associated with positive line bundles on a  projective manifold
to the case of psedoeffective singular hermitian line bundles. 

\subsection{Local measure associated with quasiplurisubharmonic functions}
\label{omega}
The content of this subsection is taken from \cite{g-z}.
Let $u,v$ be  bounded plurisubharmonic functions on some domain $D$
 in $\mathbb{C}^{n}$.
Then 
\begin{equation}\label{eq}
\mbox{\bf 1}_{\{ u > v\}}[dd^{c}u]^{n}
= \mbox{\bf 1}_{\{ u > v\}}[dd^{c}\max (u,v)]^{n}
\end{equation} 
holds in weak sense of measure in $D$ (\cite{b-t}), 
where 
\[
d^{c} := \frac{\sqrt{-1}}{4\pi}(\bar{\partial}-\partial)
\]
and 
the operation $[dd^{c}\cdot ]^{n}$ is the wedge product of 
closed positive $(1,1)$ current with bounded potential defined 
in \cite{b-t}.
Let $(X,\omega )$ be a compact K\"{a}hler manifold and we set 
\[
PSH(X,\omega ) = \{
\varphi \in L^{1}_{loc}(X)\mid \mbox{uppersemicontinuous function on $X$ such that }. 
\]
\[
\hspace{20mm} \mbox{$\omega_{\varphi}:= \omega + dd^{c}\varphi$
is a closed semipositive $(1,1)$ current on $X$}\}. 
\]
We call $PSH(X,\omega )$ the set of $\omega$ plurisubharmonic functions. 
Let $\varphi\in PSH(X,\omega )$ be a $\omega$ plurisubharmonic function on 
$X$.  

The purpose of this subsection is to define the Monge-Amp\`{e}re measure
associated with the closed positive current $\omega + dd^{c}\varphi$. 

We set 
$\varphi_{j}:= \max (\varphi , -j) \in PSH(X,\omega )$. We call the sequence
$\{ \varphi_{j}\}$ the canonical approximation of $\varphi$ by bounded
$\omega$-plurisubharmonic functions. 
This is a decreasing sequence and by (\ref{eq})
\[
\mbox{\bf 1}_{\{\varphi_{j}> k\}}[\omega + dd^{c}\varphi_{j}]^{n}
= \mbox{\bf 1}_{\{ \varphi >-k\}}[\omega + dd^{c}\max (\varphi_{j},-k)]^{n}
\]
holds. 
If $j > k$ holds, then $\{ \varphi_{j} > -k\} = \{ \varphi > -k\}$ 
and $\max (\varphi_{j},-k) = \varphi_{k}$ hold.  Hence
\[
j\geqq k \Rightarrow \mbox{\bf 1}_{\{\varphi > -j\}}
[\omega + dd^{c}\varphi_{j}]^{n}
\geqq 
\mbox{\bf 1}_{\{\varphi > -k\}}
[\omega + dd^{c}\varphi_{k}]^{n}
\]
holds. 
Since $\{\varphi > -k\} \subseteq \{ \varphi > -j\}$ holds,
we have 
\[
j\geqq k \Rightarrow \mbox{\bf 1}_{\{\varphi > -j\}}
[\omega + dd^{c}\varphi_{j}]^{n}
\geqq 
\mbox{\bf 1}_{\{\varphi > -k\}}
[\omega + dd^{c}\varphi_{k}]^{n}.
\]
We set 
\[
d\mu_{\varphi} := \lim_{j\rightarrow\infty}\mbox{\bf 1}_{\{\varphi > -j\}}
[\omega + dd^{c}\varphi_{j}]^{n}.
\]
This is a positive Borel measure which is precisely the non-pluripolar part of 
$(\omega + dd^{c}\varphi )^{n}$. 

In general the total mass
\[
\int_{X}d\mu_{\varphi}
\]
can take any value in $[0,\int_{X}\omega^{n}]$. 
This phenomena is caused by the escape of the measure toward the pluripolar 
set of $\varphi$.  

\subsection{Local volume of pseudoeffective singular hermitian line bundles}

Let $X$ be a projective manifold of dimension $n$  and let $(L,h_{L})$ be a 
pseudoeffective singular hermitian line bundle on $X$.  

\begin{definition}\label{volume}
In the above notation, we define the number $\mu (L,h_{L})$ by 
\[
\mu (L,h_{L}):= n!\limsup_{m\rightarrow\infty}m^{-n}h^{0}(X,{\cal O}_{X}(mL)
\otimes {\cal I}(h_{L}^{m}))
\]
is called the volume of $(L,h_{L})$. 
 $\square$ 
\end{definition}
\begin{definition}
Let $(L,h_{L})$ be a  pseudoeffective signgular hermitian line bundle 
on a projective manifold $X$. 
$(L,h_{L})$ is said to be big, if the volume 
$\mu (L,h_{L})$ is positive. $\square$ 
\end{definition}

In this subsection we shall define the local version of $\mu (L,h_{L})$. 

Let $h_{0}$ be a $C^{\infty}$ hermitian metric on $L$ and let 
$\varphi \in L^{1}(X)$ be the weight function  of $h_{L}$ with 
respect to $h_{0}$,i.e., $\varphi$ is a function such that 
\[
h_{L} = e^{-\varphi}\cdot h_{0}
\]
holds. 
Let $A$ be an ample line bundle on $X$ such that 
$A - L$ is very ample.  Let $\sigma$ be a nontrivial global holomorphic section  of $A-L$ and let $h_{A-L}$ be a $C^{\infty}$ hermitian metric on $A -L$ 
such that $h_{A}:= h_{0}\cdot h_{A-L}$ is a $C^{\infty}$ hermitian metric 
on $h_{A}$ with strictly positive curvature on $X$. 
Then  
\[
h_{L}\cdot \frac{1}{\mid\sigma\mid^{2}} = (h_{0}\cdot h_{A-L})\cdot e^{-(\varphi + \log h_{A-L}(\sigma ,\sigma))} 
\]
is a singular hermitian metric on $A$. 
We set 
\[
\psi := \varphi + \log h_{A-L}(\sigma ,\sigma)
\]
and 
\[
\omega = \frac{1}{2\pi}\Theta_{h_{A}}.
\]
Then 
\[
\omega_{\psi} = \omega + dd^{c}\psi
\]
is a closed positive current on $X$ and
\[
\omega_{\psi} = \frac{1}{2\pi}\Theta_{h_{L}} + (\sigma ),
\]
where $(\sigma )$ denotes the current of integration over the divisor of 
$\sigma$. 
In particular, $\psi$ is a $\omega$-plurisubharmonic function on $X$. 
Then by the result in Section \ref{omega}, we may define 
\[
d\mu (L,h_{L}) := d\mu_{\psi}. 
\]
It is easy to see that $d\mu (L,h_{L})$ is independent of the choice of 
$A,h_{A}$, etc.  

\begin{definition}\label{local volume}
Let $(L,h_{L})$ be a pseudoeffective singular hermitian line bundle 
on a projective manifold $X$. 
We call the measure $d\mu (L,h_{L})$ the local volume of $(L,h_{L})$. $\square$
\end{definition}

\subsection{Asymptotic expansion of Bergman kernels}\label{as} 
Let $X$ be a smooth projective variety and let $(L,h_{L})$ be 
a pseudoeffective singular hermitian line bundle on $X$.
 
Let $A$ be a sufficiently ample line bundle on $X$ so that 
\[
{\cal O}_{X}(K_{X}+A + mL) \otimes {\cal I}(h^{m}_{L})
\]
is globally generated for every $m \geqq 1$.
Let $h_{A}$ be a $C^{\infty}$ hermitian metric on $A$ with strictly positive 
curvature.   
Let 
\[
K(X,K_{X}+ A + mL,h_{A}\cdot h_{L}^{m})
\]
be the Bergman kernel of $K_{X} + A +mL$ with respect to the inner product 
\[
(\sigma ,\sigma^{\prime}) := (\sqrt{-1})^{n^{2}}\int_{X}\sigma \wedge \bar{\sigma}^{\prime}\cdot h_{A}\cdot h_{L}^{m}. 
\]
The reason why we need $A$ here is that to kill the higher cohomology of 
\\ ${\cal O}_{X}(K_{X}+A+mL)\otimes{\cal I}(h_{L}^{m})$ and to localize 
the estimate below.  

Now we are interested in the asymptotics of the volume form 
\[
h_{A}\cdot h_{L}^{m}\cdot K(X,K_{X}+ A + mL,h_{A}\cdot h_{L}^{m})
\]
as $m$ tends to infinity.

We may also consider the local version of the above Bergman kernel. 
To consider the local version and to localize the estimate is quite 
crucial here.

Let $x$ be a point on $X$ and let $U$ be a coordinate neighbourhood of $x$
such that $U$ is biholomorphic to a ball. 
Let 
\[
K(U,K_{X}+ A + mL,h_{A}\cdot h_{L}^{m})
\]
be the Bergman kernel of $K_{X} + A +mL\mid_{U}$ with respect to the 
inner product 
\[
(\sigma ,\sigma^{\prime}) := (\sqrt{-1})^{n^{2}}\int_{U}\sigma \wedge \bar{\sigma}^{\prime}\cdot h_{A}\cdot h_{L}^{m}. 
\]
Then it is obvious that 
\[
K(U,K_{X}+ A + mL,h_{A}\cdot h_{L}^{m})
\geqq K(X,K_{X}+ A + mL,h_{A}\cdot h_{L}^{m})
\] 
On the other hand if we  replace $A$ by its high multiple, if necessary, 
we see tat there exists a positive constant $C$ such that 
\[
K(U,K_{X}+ A + mL,h_{A}\cdot h_{L}^{m})(x)
\leqq C\cdot h_{A}\cdot K(X,K_{X}+ 2A + mL,h_{A}^{2}\cdot h_{L}^{m})(x)
\]
holds for every  $x \in B(O,1/2)$  and every positive integer $m$. 
This estimate immediately follows from the estremal propertiy of 
the Bergman kernels, i.e.,
\[
K(U,K_{X}+ A + mL,h_{A}\cdot h_{L}^{m})(x)
= \sup \{\mid\sigma\mid^{2}(x)\mid (\sqrt{-1})^{n^{2}}\int_{U}\sigma\wedge\bar{\sigma}\cdot h_{A}\cdot h_{L}^{m} = 1\}
\] 
and similar equality for $K(X,K_{X}+ 2A + mL,h_{A}^{2}\cdot h_{L}^{m})(x)$ and H\"{o}rmander's $L^{2}$-estimates for $\bar{\partial}$ operators. 

In this way,  thanks to the presense of $A$, we may localize the estimate of  the asymptotics of the global Bergman kernels in terms of that of local Bergman kernels.    

The local asymptotics of $h_{A}\cdot h_{L}^{m}\cdot
K(U,K_{X}+ A + mL,h_{A}\cdot h_{L}^{m})$ (hence also the global asymptotics) can be explored by the 
following well known theorem, if $h_{L}$ is $C^{\infty}$. 

\begin{theorem}\label{local asymptotics}(\cite{c,ti,ze})
Let $\Omega$ be a pseudoconvex domain in $\mathbb{C}^{n}$ and 
let $\varphi$ be a $C^{\infty}$ plurisubharmonic function on $\Omega$. 
For a positive integer $m$, let $K(K_{\Omega},e^{-m\varphi})$ be the Bergman kernel  of $K_{\Omega}$ with respect to the inner product
\[
(f,g) := (\sqrt{-1})^{n^{2}}\int_{\Omega}e^{-m\varphi}\cdot f\wedge \bar{g}. 
\]
Then  
\[
e^{-m\varphi}\cdot K(K_{\Omega},e^{-m\varphi}) = \frac{(dd^{c}\varphi )^{n}}{n!}m^{n} + O(m^{n-1}) 
\] 
holds. $\square$. 
\end{theorem}
Theorem \ref{local asymptotics} can be viewed as a  (weak version of ) local Riemann-Roch theorem. 

\begin{definition}\label{local volume asymptotics}
Let $(L,h_{L})$ be a pseudoeffective singular hermitian line bundle on a projective manifold $X$ of dimension $n$. Let $A$ be a sufficiently ample line bundle
on $X$. 
We set 
\[
d\mu^{+}(L,h_{L}):= n!\cdot\limsup_{m\rightarrow\infty}m^{-n}\cdot h_{A}\cdot h_{L}^{m}\cdot K(X,K_{X} +A +mL,h_{A}\cdot h_{L}^{m})(z).
\]
$d\mu^{+}(L,h_{L})$ is said to be the  upper local volume of $(L,h_{L})$.
Similarly we set 
\[
d\mu^{-}(L,h_{L}):=  n!\cdot\liminf_{m\rightarrow\infty}m^{-n}\cdot h_{A}\cdot h_{L}^{m}\cdot K(X,K_{X} +A +mL,h_{A}\cdot h_{L}^{m})(z).
\]
$d\mu^{-}(L,h_{L})$ is said to be the lower local volume of $(L,h_{L})$.
$\square$
\end{definition}
Since 
\[
\int_{X}h_{A}\cdot h_{L}^{m}\cdot K(X,K_{X}+ A + mL,h_{A}\cdot h_{L}^{m})
= \dim H^{0}(X,{\cal O}_{X}(K_{X} + A + mL)\otimes{\cal I}(h_{A}\cdot h_{L}^{m}))
\]
holds, by Lebesgue-Fatou's lemma, we have that 
\[
\int_{X}d\mu^{-}(L,h_{L}) \leqq \mu (L,h_{L}) \leqq \int_{X}d\mu^{+}(L,h_{L})
\]
hold. 
\subsection{A general conjecture for the asymptotic expansion of Bergman kernels}
Let $X$ be a smooth projective variety of dimension $n$ and let $(L,h_{L})$ be a pseudoeffective singular hermitian line bundle on $X$. 
Let $A$ be a sufficiently ample line bundle on $X$ and let $h_{A}$ be 
a $C^{\infty}$ hermitian metric on $A$ with strictly positive curvature. 
\begin{conjecture}\label{gen volume}
\begin{enumerate}
\item $d\mu^{+}(L,h_{L}) = d\mu^{-}(L,h_{L})$ holds (cf. Definition \ref{local volume asymptotics}).
In particular \\
$n!\cdot\lim_{m\rightarrow\infty}m^{-n}\cdot h_{A}\cdot h_{L}^{m}\cdot K(X,K_{X}+A+mL,h_{A}\cdot h_{L}^{m})$ exists on $X$. 
\item
\[
n!\cdot\lim_{m\rightarrow\infty}m^{-n}\cdot h_{A}\cdot h_{L}^{m}\cdot K(X,K_{X}+A+mL,h_{A}\cdot h_{L}^{m}) = d\mu (L,h_{L}), 
\]
where $d\mu (L,h_{L})$ denotes the local volume defined in Definition \ref{local volume}.
\end{enumerate} 
$\square$  
\end{conjecture}
Conjecture \ref{gen volume} is true for the case that $h_{L}$ is smooth
by \cite{ti,ze}.   Also by essentially the same proof, it is easy to verify that Conjecture \ref{gen volume} is true for the case that $h_{L}$ has algebraic singularities.  
But in general, Conjecture \ref{gen volume} seems to be very difficult to prove, since the asymptotic expansion of Bergman kernels breaks down, if 
the hermitian metric is not $C^{\infty}$.
Of course one may approximate $h_{L}$ by a singular hermitian metric 
with algebraic singularities (\cite{dem}).  But the difficulty arises when we take the 
limit, i.e., the lower order term of the asymptotic expansion may blow up. 
 
In the next subsection, we shall prove a lower estimate for the asymptotic 
expansion of Bergman kernels, when the curvature is strictly positive.    

\subsection{Asymptotic expansion of Bergman kernels of singular hermitian line bundles 
with strictly positive curvature currents}

In this subsection, we shall estimate the asymptotic expansion of Bergman kernels of singular hermitian line bundles with strictly positive curvature.
We prove that the coefficient of the top term of the expansion at a point  is strictly 
positive when the curvature is bounded from below by a strictly positive form locally 
around the point. 
\begin{theorem}\label{MA}
Let $\Omega$ be a bounded pseudoconvex domain in $\mathbb{C}^{n}$ 
contained in the unit open ball $B(O,1)$ centered at the origin with radius $1$   and 
let $\varphi$ be a plurisubharmonic function on $\Omega$ such that   
\[
dd^{c}\varphi \geqq - dd^{c}\log (1 - \parallel z\parallel^{2}) 
\]
holds on $\Omega$, where $\parallel z\parallel^{2} = \sum_{i=1}^{n}\mid z_{i}\mid^{2}$. 
Then there exists a positive constant $C$ independent of 
$\varphi$ such that  for every positive integer $m$ 
\[
e^{-m\varphi}\cdot K(\Omega,e^{-m\varphi}) \geqq C\cdot m^{n}\cdot \mid dz_{1}\wedge\cdots\wedge dz_{n}\mid^{2}
\]
holds on $\Omega$. $\square$
\end{theorem}
{\bf Proof}. 
Since $B(O,1)$ is homogeneous, we may assume that $\Omega$ contains the origin
$O\in \mathbb{C}^{n}$ and it is is enough to prove that 
\[
e^{-m\varphi}\cdot K(\Omega,e^{-m\varphi})(O) \geqq C\cdot m^{n}\cdot 
\mid dz_{1}\wedge\cdots\wedge dz_{n}\mid^{2}
\]
holds at the origin $O$.  
Let $\rho \in C^{\infty}_{0}(\mathbb{C}^{n})$ be a function such that 
\begin{enumerate}
\item $0\leqq \rho \leqq 1$ hold on $\mathbb{C}^{n}$.
\item $\mbox{Supp}\,\rho \subset B(O,1/2)$ holds.  
\item $\mid d\rho\mid < 3$, where the norm is taken with respect to 
the standard Euclidean metric on $\mathbb{C}^{n}$.
\end{enumerate}
Since $dd^{c}\log (1 - \parallel z\parallel^{2})$ is a K\"{a}hler form 
invariant under $\mbox{Aut}(B(O,1))$, for every sufficiently large $m$,   
there exists a positive constant $c > 1$ depending only on $\rho$ such that 
\begin{equation}\label{cut}
m\cdot dd^{c}\log (1 - \parallel z\parallel^{2}) 
+ dd^{c}(\rho (\frac{c}{\sqrt{m}}z)\cdot\frac{1}{\parallel z\parallel^{2n}}) 
\geqq -\frac{m}{2}\cdot dd^{c}\log (1 - \parallel z\parallel^{2}) 
\end{equation}
holds on  $B(O,1)$. 
Hence by the standard $L^{2}$-estimates, we have the following lemma.
\begin{lemma}\label{local}
There exists a positive constant $c_{0}$ such that for every sufficiently large $m \geqq 1$,  
\[
K(\Omega ,e^{-m\varphi})(O) \geqq c_{0}\cdot K(B(O,\frac{c}{\sqrt{m}}),e^{-m\varphi})(O)
\]
holds.    $\square$
\end{lemma}
{\bf Proof of Lemma \ref{local}.}
Let $f$ be a holomorphic $(n,0)$ form on $B(O,\frac{c}{\sqrt{m}})$ such that 
\[
\mid f\mid^{2}(O) = K(B(O,\frac{c}{\sqrt{m}}),e^{-m\varphi})(O)
\]
and 
\[
\int_{B(O,\frac{c}{\sqrt{m}})}\mid f\mid^{2}\cdot e^{-m\varphi}d\mu = 1
\]
holds, i.e., $f$ is a peak section at $O$ with respect to the metric 
$e^{-m\varphi}$. 

Let us fix a complete K\"{a}hler metric $g$ on $\Omega$ and let $\omega_{g}$ denote the K\"{a}hler form associated with $g$.
Let $\gamma(z)$ be the minimum eigenvalue of  $dd^{c}\log (1 - \parallel z\parallel^{2})$ with respect to $g(z)$.   Then 
\[
dd^{c}\log (1 - \parallel z\parallel^{2}) 
\geqq  \gamma(z)\omega_{g}(z)
\] 
holds. 
Then by the $L^{2}$-estimates and (\ref{cut}), we have that
there exists a $(n,0)$ form $u$ on $\Omega$ such that  
\[
\bar{\partial}u = \bar{\partial}(\rho(\frac{c}{\sqrt{m}}z)f(z))
\]
and 
\begin{equation}\label{estimate}
\int_{\Omega}\mid u\mid^{2}\cdot e^{-m\varphi}
\cdot e^{\rho(\frac{c}{\sqrt{m}}z)}\frac{1}{\parallel z\parallel^{2n}}
\leqq \frac{2}{m}\int_{\Omega}\frac{1}{\gamma(z)}
\mid \bar{\partial}(\rho(\frac{c}{\sqrt{m}}z)f(z))\mid^{2}\cdot e^{\rho(\frac{c}{\sqrt{m}}z)}\cdot \frac{1}{\parallel z\parallel^{2n}}d\mu_{g}
\end{equation}
hold, where $\mid \bar{\partial}(\rho(\frac{c}{\sqrt{m}}z)f(z))\mid$ denotes 
the norm with respect to $g$ and $d\mu_{g}$ denotes the volume form with respect to $g$.  
Since there exists a positive constant $C_{1}$ independent of $m$ such that  
\[
\int_{\Omega}\frac{1}{\gamma(z)}
\mid \bar{\partial}(\rho(\frac{c}{\sqrt{m}}z)f(z))\mid^{2}\cdot e^{\rho(\frac{c}{\sqrt{m}}z)}\cdot \frac{1}{\parallel z\parallel^{2n}}d\mu_{g}
\leqq C_{1}\cdot m
\]
holds for every sufficiently large $m$, by (\ref{estimate}), we see that 
there exists a positive constant $C_{2}$ independent of $m$ such that 
\[
\int_{\Omega}\mid u\mid^{2}\cdot e^{-m\varphi}
\cdot e^{\rho(\frac{c}{\sqrt{m}}z)}\frac{1}{\parallel z\parallel^{2n}}
\leqq C_{2}. 
\]
Since $u(O) = 0$ holds by the construction, we see that 
\[
\rho(\frac{c}{\sqrt{m}}z)f(z) - u
\]
is a holomorphic extension of $f(O)$ at $O$ and there exists 
a positive constant $C_{3}$ independent of $m$ such that  
\[
\int_{\Omega}\mid\rho(\frac{c}{\sqrt{m}}z)f(z) - u\mid^{2}\cdot e^{-m\varphi}
\leqq C_{3}
\]
holds.  Hence by the extremal property of Bergman kernels, we see that 
\[
K(\Omega ,e^{-m\varphi})(O) \geqq \frac{1}{C_{3}}\cdot K(B(O,\frac{c}{\sqrt{m}}),e^{-m\varphi})(O)
\]
holds.  This completes the proof of Lemma \ref{local}.  
$\square$  \vspace{3mm} \\
On the other hand, by the $L^{2}$-extension theorem,  we see that
there exists a positive constant $C_{0}$ independent of $m$ such that 
for every sufficiently large $m$,  
\[
K(B(O,\frac{c}{\sqrt{m}}),e^{-m\varphi})(O) \geqq C_{0}\cdot e^{m\varphi(O)}\cdot m^{n}
\cdot \mid dz_{1}\wedge\cdots\wedge dz_{n}\mid^{2}
\]
holds. 
Hence combining the above inequality and Lemma \ref{local}, we see that 
\[
K(\Omega ,e^{-m\varphi})(O) \geqq c_{0}\cdot C_{0}\cdot e^{m\varphi(O)}\cdot m^{n}
\cdot \mid dz_{1}\wedge\cdots\wedge dz_{n}\mid^{2}
\]
holds.   This completes the proof of Theorem \ref{MA}. \vspace{3mm} $\square$ 

\noindent The global version of Theorme \ref{MA} is as follows. 
\begin{theorem}\label{global}
Let $X$ be a smooth projective variety of dimension $n$ and let $(L,h_{L})$ be a singular 
hermitian line bundle on $X$.
Let $\omega$ be a $C^{\infty}$ K\"{a}hler form on $X$.    
Suppose that there exists a positive constant $\varepsilon$ such that 
\[
\Theta_{h_{L}} \geqq \varepsilon\cdot\omega 
\]
holds on $X$.   
Then there exists a positive constant $C$ such that  for every sufficiently large positive integer $m$, 
\[
h_{L}^{m}\cdot K(X,K_{X}+mL,h_{L}^{m}) \geqq  C\cdot m^{n}\cdot \omega^{n}
\]
holds on $X$.    $\square$ 
\end{theorem}
The proof of Theorem \ref{global} follows from the local version (Theorem \ref{MA}) by the localization principle (cf. Section \ref{as}). 

\section{Kodaira's lemma for big pseudoeffective singular hermitian line bundles}
In this section we shall prove an analogue of Kodaira's lemma 
for big pseudoeffecive singular hermitian line bundles.
Kodaira's lemma has been extensively used in algebraic geometry
(cf. \cite{k-o,ka1}). 
To prove Theorem \ref{main}, we need an analogue of 
Kodaira's lemma for big pseudoeffective singular hermitian line bundles.

Although I can state the new version as a lemma, it will be useful to state the singular hermitian version  as a theorem.  Because I believe that 
the new version will be also fundamental in complex geometry.

\subsection{Statement of the theorem}
First we shall state the original Kodaira's lemma. 
\begin{theorem}(\cite[Appendix]{k-o},\cite[Lemma]{ka1})\label{original}
Let $X$ be a smooth projective variety and let $D$ be a big divisor on $X$.
Then there exists an effective $\mathbb{Q}$-divisor $E$ such that 
$D - E$ is an ample $\mathbb{Q}$-divisor. $\square$
\end{theorem}
The analogue for the case of big pseudoeffective singular hermitian line bundles is  stated as follows. 

\begin{theorem}\label{kodaira}
Let $X$ be a projective manifold and let $(L,h_{L})$ be a big 
psedoeffective singular hermitian line bundle. 
Then there exists a singular hermitian metric $h^{+}_{L}$ on $L$ such that 
\begin{enumerate}
\item $\Theta_{h^{+}_{L}}$ is strictly positive 
everywhere on $X$, 
\item $h^{+}_{L} \geqq h_{L}$ holds on $X$. 
\end{enumerate} $\square$
\end{theorem}
Let us explain the relation between Theorems  \ref{original} and \ref{kodaira}. 
\noindent Let $D,E$ be as in Theorem \ref{original}.  Let us identify divisors with line bundles. 
Theorem \ref{original} says that there exists a $C^{\infty}$ hermitian metrics 
$h_{D},h_{E}$ on $D,E$ respectively (the notion of hermitian metrics 
naturally extends to the case of $\mathbb{Q}$-line bundles) such that 
the curvature of $h_{D}\cdot h_{E}^{-1}$ is stricly positive. 
Let $\sigma_{E}$ be a multivalued holomorphic section of $E$ with divisor $E$
such that $h_{E}(\sigma_{E},\sigma_{E})\leqq 1$ on $X$. 
Then 
\[
h_{D}^{+}:= \frac{h_{D}}{h_{E}(\sigma_{E},\sigma_{E})}
\]
is a singular hermitian metric on $D$  such that 
\begin{enumerate}
\item $\Theta_{h_{D}^{+}}$
is strictly positive everywhere on $X$.
\item $h_{D} \leqq h_{D}^{+}$
holds on $X$.
\end{enumerate}
In this way Theorem \ref{kodaira} can be viewed as an analogue of the usual 
Kodaira's lemma to the case of big pseudoeffective singular hermitian line bundles.
\subsection{Proof of Theorem \ref{kodaira}}
The proof of Theorem \ref{kodaira} presented here is not very much different from the original proof of Kodaira's lemma (cf. \cite{ka1} or \cite[Appendix]{k-o}).  But it requires estimates of  Bergman kernels  and  additional care for the multiplier 
ideal sheaves. 

Let $X$ be a smooth projective variety of dimension $n$ and let 
$(L,h_{L})$ be a big pseudoeffective singular hermitian line bundle on $X$. 
Let $\omega$ be a K\"{a}hler form on $X$ and let $dV$ be the associated 
volume form on $X$. 
Let $H$ be a smooth very ample divisor on $X$.  
The following lemma is a singular hermitian version of the theorem in \cite{tr}. 
\begin{lemma}\label{big}
There exists a positive integer $m_{0}$ such that 
$m_{0}(L,h_{L}) - H$ is big, i.e.,  
\[
\limsup_{\ell\rightarrow\infty}\ell^{-n}\cdot\dim H^{0}(X,{\cal O}_{X}(\ell (m_{0}L - H)\otimes {\cal I}(h_{L}^{m_{0}\ell})) > 0
\]
holds.  $\square$ 
\end{lemma}
{\bf Proof of Lemma \ref{big}.}
Replacing $H$ by a suitable member of $\mid\! H\!\mid$, by Lemma \ref{generic},
we may assume that 
\[
{\cal I}(h_{L}^{m})\mid_{H} = {\cal I}(h_{L}^{m}\mid_{H})
\]
holds for every $m \geqq 1$. 
Let us consider the exact sequence 
\[
0\rightarrow H^{0}(X,{\cal O}_{X}(mL - H)\otimes {\cal I}(h_{L}^{m}))\rightarrow H^{0}(X,{\cal O}_{X}(mL)\otimes {\cal I}(h_{L}^{m})) 
\]
\[
\hspace{45mm}\rightarrow H^{0}(H,{\cal O}_{H}(mL)\otimes {\cal I}(h_{L}^{m}\mid_{H})). 
\]
Then since $\mu (L,h_{L}) > 0$ and 
\[
\dim H^{0}(H,{\cal O}_{H}(mL)\otimes {\cal I}(h_{L}^{m}\mid_{H})) = O(m^{n-1})
\]
we see that for every sufficiently large 
$m$, 
\[
H^{0}(X,{\cal O}_{X}(mL - H)\otimes {\cal I}(h_{L}^{m}))\neq 0
\]
holds. 

To prove Lemma \ref{big}, we need to refine the above argument a little bit. 
Let $m_{0}$ be a positive integer such that 
\begin{equation}\label{m}
m_{0} >n\cdot \frac{(L,h_{L})^{n-1}\!\!\cdot H}{(L,h_{L})^{n}} 
\end{equation}
holds.   For very general $H_{1}^{(\ell)},\cdots H^{(\ell)}_{\ell}\in \mid\! H\!\mid$, by Lemma \ref{generic}, replacing $m$ by $m_{0}\ell$ and 
$H$ by $\ell H$, we have the exact sequence 
\[
0\rightarrow H^{0}(X,{\cal O}_{X}(\ell (m_{0}L - H))\otimes{\cal I}(h_{L}^{m_{0}\ell}))\rightarrow H^{0}(X,{\cal O}_{X}(m_{0}\ell L)\otimes {\cal I}(h_{L}^{m_{0}\ell})).  
\]
\[
\hspace{40mm} \rightarrow 
\oplus_{i=1}^{\ell}H^{0}(H^{(\ell)}_{i},{\cal O}_{H_{i}}(m_{0}\ell L)\otimes {\cal I}(h_{L}^{m_{0}\ell}\mid_{H_{i}})). 
\]
We note that $\{ H_{i}^{(\ell)}\}_{i=1}^{\ell}$ are  chosen for each $\ell$.
If we take $\{ H_{i}^{(\ell)}\}_{i=1}^{\ell}$ very general, 
we may assume  that 
\[
\dim H^{0}(H^{(\ell)}_{i},{\cal O}_{H_{i}}(mL)\otimes {\cal I}(h_{L}^{m}\mid_{H_{i}}))
\]
is independent of $1\leqq i\leqq \ell$ for every $m$. 
This implies that 
\[
\limsup_{\ell\rightarrow\infty}\ell^{-n}\cdot\dim H^{0}(X,{\cal O}_{X}(\ell (m_{0}L - H))\otimes {\cal I}(h_{L}^{m_{0}\ell}))
\]
\[
\hspace{20mm}
\geqq \frac{1}{n!}(L,h_{L})^{n}\cdot m_{0}^{n}- \frac{1}{(n-1)!}
 \{(L,h_{L})^{n-1}\!\!\cdot H\} \cdot m_{0}^{n-1}
\]
holds.  By (\ref{m}), we see that 
\[
\frac{1}{n!}(L,h_{L})^{n}\cdot m_{0}^{n}- \frac{1}{(n-1)!}\{(L,h_{L})^{n-1}\!\!\cdot H\}m_{0}^{n-1}
\]
is positive. 
This completes the proof of Lemma \ref{big}. $\square$
\vspace{5mm}\\
Let $A$ be a sufficiently ample line bundle on $X$ and let $h_{A}$ be a 
$C^{\infty}$ hermitian metric such that the curvature of $h_{A}$ is everywhere
strictly positive on $X$.  Here the meaning of ``sufficiently ample'' will 
be specified later. 
Let $m$ be a positive integer. 
Let us consider the inner product 
\[
(\sigma ,\sigma^{\prime}) := \int_{X}h_{A}\cdot h_{L}^{m}\cdot \sigma
\cdot\bar{\sigma}^{\prime}\, dV
\]
on $H^{0}(X,{\cal O}_{X}(A+ mL)\otimes {\cal I}(h_{L}^{m}))$ and 
let $K_{m}$ be the associated (diagonal part of) Bergman kernel.
Let us consider the subspace 
\[ 
H^{0}(X,{\cal O}_{X}(A+\ell (m_{0}L - H))\otimes {\cal I}(h_{L}^{m_{0}\ell}))
\subset H^{0}(X,{\cal O}_{X}(A+m_{0}\ell L)\otimes {\cal I}(h_{L}^{m_{0}\ell}))
\]
as a Hilbert subspace and let $K_{m_{0}\ell}^{+}$ denotes the associated 
Bergman kernel with respect to the restriction of the inner product 
on \\ $H^{0}(X,{\cal O}_{X}(A + m_{0}\ell L)\otimes {\cal I}(h_{L}^{m_{0}\ell}))$
to the subspace
$H^{0}(X,{\cal O}_{X}(A+\ell (m_{0}L - H))\otimes {\cal I}(h_{L}^{m_{0}\ell}))
$. 
Then by definition, we have the trivial inequality :
\begin{equation}\label{trivial}
K_{m_{0}\ell}^{+} \leqq K_{m_{0}\ell}
\end{equation}
holds on $X$ for every $\ell \geqq 1$. 

The next lemma follows from the same argument as in \cite{dem}.  

\begin{lemma}\label{reciprocity}(\cite{dem})
If $A$ is sufficiently ample,
\[
h_{L} := \mbox{the lower envelope of}\,\,(\limsup_{m\rightarrow\infty}\sqrt[m]{K_{m}})^{-1}.
\]
holds. $\square$
\end{lemma}
{\bf Proof of Lemma \ref{reciprocity}}.
The proof has been given in \cite{dem}.  But for the completeness, 
we shall reproduce the proof here. 

By the $L^{2}$-extension theorem  (Theorem \ref{o-t} or Theorem \ref{extension}), if $A$ is sufficiently ample, 
there exists a positive constant $C_{0}$ such that  
\begin{equation}\label{l}
K_{m} \geqq C_{0}\cdot h_{A}^{-1}\cdot h_{L}^{-m}
\end{equation}
holds on $X$ for every $m$.

On the other hand,  let $x\in X$ be an arbitrary point and let 
$(U,z_{1},\cdot ,z_{n})$ be a coordinate neighbourhood centered at $x$ 
such that $U$ is biholomorphic to the open unit ball $B(O,1)$ in $\mathbb{C}^{n}$ centered at the origin via the coordinate.  
Taking $U$ to be sufficiently small, we may and do assume that $\mbox{\bf e}_{A},\mathbb{e}_{L}$ be the holomorphic 
frame of $A$ and $L$ on $U$ respectively. 
Then with respect to these frame, we may express $h_{A},h_{L}$ as 
\[
h_{A} = e^{-\varphi_{A}}, h_{L} = e^{-\varphi_{L}}
\]
respectively in terms of plurisubharmonic functions $\varphi_{A},\varphi_{L}$
on $U$. 
By the extremal property of Bergman kernels, we see that 
\[
K_{m}(x) = \sup \{ \mid\sigma (x)\mid^{2}\mid 
\sigma \in \Gamma (X,{\cal O}_{X}(A + mL)), \int_{X}\mid\sigma\mid^{2}\cdot h_{A}\cdot h_{L}^{m}\cdot dV = 1\}.
\]
Let $\sigma_{0} \in \Gamma (X,{\cal O}_{X}(A + mL))$ with 
\[
\int_{X}\mid\sigma_{0}\mid^{2}\cdot h_{A}\cdot h_{L}^{m}\cdot dV = 1
\]
and $\mid\sigma_{0}(x)\mid^{2} = K_{m}(x)$.
Let us write $\sigma_{0} = f\cdot\mbox{\bf e}_{A}\cdot\mbox{\bf e}_{L}^{m}$
on $U$ by using a holomorphic function $f$ on $U$.
By the submeanvalue property of plurisubharmonic functions, we have that 
\begin{eqnarray*}
\mid f(O)\mid^{2} &\leqq &  \frac{1}{\mbox{vol}(B(O,\varepsilon ))}\int_{B(O,\varepsilon)}\mid f\mid^{2}d\mu \\
&\leqq & (\sup_{B(O,\varepsilon )}e^{\varphi_{A}}\cdot e^{m\varphi_{L}})
\cdot (\frac{1}{\mbox{vol}(B(O,\varepsilon ))}\int_{B(,\varepsilon )}\mid f\mid^{2}e^{-\varphi_{A}}\cdot e^{-m\varphi_{L}}dV)\cdot 
(\sup_{B(O,\varepsilon )}\frac{d\mu}{dV})  \\
&\leqq & \frac{1}{\mbox{vol}(B(O,\varepsilon ))}\cdot (\sup_{B(O,\varepsilon )}e^{\varphi_{A}}\cdot e^{m\varphi_{L}})
\cdot (\sup_{B(O,\varepsilon )}\frac{d\mu}{dV})
\end{eqnarray*}
hold, where $d\mu$ is the standard Lebesgue measure on $\mathbb{C}^{n}$. 
Hence there exists a positive constant $C_{\varepsilon}$ independent of $m$
\begin{equation}\label{up}
K_{m}(x) \leqq C_{\varepsilon}\cdot \sup_{w\in B(O,\varepsilon )}(h_{A}^{-1}\cdot h_{L}^{-m})(w)\cdot dV
\end{equation}
holds.  
By (\ref{l}) and (\ref{up}), we see that 
\[
C^{\frac{1}{m}}(h_{A}^{-\frac{1}{m}}h_{L}^{-1}) \leqq  \sqrt[m]{K_{m}(x)} \leqq C_{\varepsilon}^{\frac{1}{m}}\cdot 
(\sup_{w\in B(O,\varepsilon )}(h_{A}^{-1}\cdot h_{L}^{-m})(w)\cdot dV)^{\frac{1}{m}}   
\]
holds. 
Hence letting $m$ tend to infinity and then letting 
$\varepsilon$ tend to $0$, we have the desired equality. 
$\square$\vspace{5mm} \\
We note that 
\[
\int_{X}h_{A}\cdot h_{L}^{m}\cdot K_{m}\cdot dV = \dim H^{0}(X,{\cal O}_{X}(A+ mL)\otimes {\cal I}(h_{L}^{m}))
\]
and 
\[
\int_{X}h_{A}\cdot h_{L}^{m_{0}\ell}\cdot K_{m_{0}\ell}^{+}\cdot dV = \dim H^{0}(X,
{\cal O}_{X}(A+\ell (m_{0}L - H))\otimes {\cal I}(h_{L}^{m_{0}\ell}))
\]
hold.  
Hence by Lemma \ref{big}
\begin{equation}\label{ve}
\limsup_{\ell\rightarrow\infty}(m_{0}\ell)^{-n}\cdot \int_{X}h_{L}^{m_{0}\ell}\cdot K^{+}_{m_{0}\ell}\cdot dV  > 0
\end{equation}
holds. 
Then by Fatou's lemma, we see that 
\[
\int_{X}\limsup_{\ell\rightarrow\infty}\frac{h_{A}\cdot h_{L}^{m_{0}\ell}\cdot K^{+}_{m_{0}\ell}}{(m_{0}\ell )^{n}}
\geqq
\limsup_{\ell\rightarrow\infty}\int\frac{h_{A}\cdot h_{L}^{m_{0}\ell}\cdot K^{+}_{m_{0}\ell}}{(m_{0}\ell)^{n}} > 0 
\] 
hold.
In particular
\[
\limsup_{\ell\rightarrow\infty}\frac{h_{A}\cdot h_{L}^{m_{0}\ell}\cdot K^{+}_{m_{0}\ell}}{(m_{0}\ell )^{n}}
\]
is not identically $0$\footnote{At this moment, there is a possibility that it is 
identically $+\infty$.}.
This implies that 
\[
\limsup_{\ell\rightarrow\infty}\sqrt[m_{0}\ell]{K^{+}_{m_{0}\ell}}
\]
is not identically $0$ and by Lemma \ref{reciprocity} and (\ref{trivial}), 
it is finite.  
\noindent Let $h_{H}$ be a $C^{\infty}$ hermitian metric on $H$ with strictly positive
curvature and let $\tau$ be a global holomorphic section of ${\cal O}_{X}(H)$
with divisor $H$ such that $h_{H}(\tau ,\tau) \leqq 1$ holds on $X$. 
We set 
\[
h_{L}^{+} := (\limsup_{\ell\rightarrow\infty}\sqrt[m_{0}\ell]{K^{+}_{m_{0}\ell}}\,\,)^{-1}\cdot h_{H}(\tau ,\tau ).
\]
Then $h_{L}^{+}$ is a singular hermitian metric
on $L$, since \\ $(\limsup_{\ell\rightarrow\infty}\sqrt[m_{0}\ell]{K^{+}_{m_{0}\ell}}\,\,)^{-1}\cdot \mid\tau\mid^{2}$ can be viewed as  a singular hermitian metric 
on $L - H$ with semipositive curvature current.  By the construction it is clear that the curvature current of $h_{L}^{+}$ is bigger than or equal to the curvature of $h_{H}$. 
In particular the curvature current of $h_{L}^{+}$ is strictly positive. 
And by the construction 
\[
h_{L} \leqq h_{L}^{+}
\]
holds on $X$. 
This completes the proof of Theorem \ref{kodaira}. $\square$

\section{Dynamical construction of an AZD}\label{Dy}
This section is almost completely the same as \cite[Section 3]{tu}. 
The only difference is that we consider the adjoint line bundles instead
of canonical bundles. 
\subsection{Sub extension problem of a singular hermitian metric with 
semipositive curvature}
Let $X$ be a smooth projective variety and let $S$ be a smooth subvariety 
of $X$.   Let $E$ be a line bundle on $X$.  
Suppose that there exists a singular hermitian metric $h_{E,S}$ on $E\mid_{S}$
such that $\Theta_{h_{E,S}} \geqq 0$ holds on $S$. 
We shall consider the following sub extension problem of singular hermitian metrics. 

\begin{problem}\label{extmetric}
Let $X,S,E,h_{E,S}$ be as above. 
Construct a singular hermitian metric $h_{E}$ on $E$ such that 
\begin{enumerate}
\item $h_{E}\mid_{S} \leqq h_{E,S}$ holds on $S$.  
\item $\Theta_{h_{E}}\geqq 0$ holds on $X$. 
\end{enumerate}
$\square$ 
\end{problem}

\noindent Of course  in general such a sub extension $h_{E}$ does not exist. 
But under certain conditions, sub extension $h_{E}$ exists. 
 Since singular hermitian metrics are real object, it is not easy to 
extend them directly.

The sub extension strategy used in this article is as follows. 
\begin{enumerate}
\item Take a sufficiently ample line bundle $A$ on $X$. 
\item Approximate $h_{E,S}$ by a sequence of singular hermitian metrics 
$\{ h_{m,S}\}$, where $h_{m,S}$ is an algebraic singular hermitian metric  (cf. Definition \ref{algsing}) 
on $\frac{1}{m}A\mid_{S} + L\mid_{S}$ which is 
of the form
\[
h_{m,S} =  (\sum_{i=0}^{N(m)}\mid \tau_{i}\mid^{2})^{-\frac{1}{m}}, 
\]
where  $\tau_{i} \in H^{0}(S,{\cal O}_{S}(A + mL))$. 
\item Find a holomorphic extension $\tilde{\tau}_{i} \in H^{0}(X,{\cal O}_{X}(A + mL))$  of $\tau_{i}$ for every $i= 0,\cdots ,N(m)$. 
\item Define  an algebraic singular hermitian metric $\tilde{h}_{m}$ 
on $\frac{1}{m}A + L$ by 
\[
\tilde{h}_{m}:= (\sum_{i=0}^{N(m)}\mid \tilde{\tau}_{i}\mid^{2})^{-\frac{1}{m}}.
\]
\item Prove the existence of 
\[
h_{E} := \mbox{the lower envelope of}\,\,\,\,\liminf_{m\rightarrow\infty}\tilde{h}_{m}
\]
as a nontrivial singular hermitian metric on $E$. 
\end{enumerate}

\noindent In the above strategy one cannot expect that the equality  
$h_{E}\mid_{S} = h_{E,S}$ holds on $S$, since we have taken the lower envelope. But in most applications, sub extension is enough.

\subsection{Dynamical construction of singular hermitian metrics}

To implement the sub extension strategy in Section 6.1, first task 
is to approximate the given singular hermitian metric by a sequence of 
algebraic singular hermitian metrics.
Here we shall consider the case that the given hermitian metric 
is an AZD. 

Instead of approximating the singular hermitian metric, 
we shall consruct another AZD which is a limit of the 
algebraic singular hermitian metrics.  
  
Let $X$ be a smooth projective variety and let 
$K_{X}$ be the canonical line bundle of $X$.
Let $n$ denote the dimension of $X$.  
Let $(L,h_{L})$ be a pseudoeffective singular hermitian line bundle on $X$.
Let $dV$ be a $C^{\infty}$ volume form on $X$.    
Suppose that $(K_{X}+ L, dV^{-1}\cdot h_{L})$ is weakly pseudoeffective (cf. Definition \ref{wpe}).  
Then 
\[
E(K_{X}+ L, dV^{-1}\cdot h_{L}):= \{ \varphi\in L^{1}_{loc}(X) \mid  \varphi \leqq 0, \,\,\Theta_{dV^{-1}\cdot h_{L}} + \sqrt{-1}\partial\bar{\partial}\varphi \geqq 0 \}
\]
is nonempty (cf. the proof of Theorem \ref{AZD2}).
Then 
$(K_{X} + L,dV^{-1}\cdot h_{L})$ admits an AZD $h$ as in Theorem \ref{AZD2}, i.e.,
the following holds :  
\begin{enumerate}
\item $\Theta_{h}\geqq 0$,  
\item 
$H^{0}(X,{\cal O}_{X}(m(K_{X} + L))\otimes {\cal I}_{\infty}(h^{m}))
\simeq H^{0}(X,{\cal O}_{X}(m(K_{X} + L))\otimes {\cal I}_{\infty}(h_{L}^{m}))
$
holds for every $m\geqq 0$. 
\end{enumerate}
We assume that for every ample line bundle $B$ on $X$
\[
\dim H^{0}(X,{\cal O}_{X}(B+ m(K_{X}+L))\otimes {\cal I}(h^{m}))
= O(m^{\nu})
\]
holds, where  $\nu$ denotes the numerical Kodaira dimension $\nu_{num}(K_{X}+ L,h)$ of $(K_{X} + L,h)$ (cf. Definition \ref{numerical}).
  We note that if $h$ is normal, this follows from Theorem \ref{nak}. 

Let $A$ be a sufficiently ample line bundle on $X$ 
such that for every pseudoeffective singular hermitian 
line bundle $(F,h_{F})$, 
\[
{\cal O}_{X}(A+F)\otimes{\cal I}(h_{F})
\]
and 
\[
{\cal O}_{X}(K_{X}+A+F)\otimes{\cal I}(h_{F})
\]
are globally generated. 
This is possible by the $L^{2}$-estimate of $\bar{\partial}$ operator (cf. \cite[p. 667, Proposition 1]{si}). 

Let $h_{A}$ be a $C^{\infty}$ hermitian metric on $A$
 with strictly positive curvature. 
Let $\{ \sigma^{(1)}_{0},\cdots ,\sigma^{(1)}_{N(1)}\}$
be a complete orthonormal basis of \\
$H^{0}(X,{\cal O}_{X}(K_{X}+L+A)\otimes{\cal I}(h_{L}))$ with 
respect to the inner product 
\[
(\sigma ,\tau ) := (\sqrt{-1})^{n}\int_{X}\sigma\wedge \bar{\tau}\cdot h_{L}\cdot h_{A}
\]
We set 
\[
K_{1} := \sum_{i=0}^{N(0)}\mid\sigma_{i}^{(1)}\mid^{2}.
\]
We define the singular hermitian metric $h_{1}$ on $K_{X}+L+A$ by  
\[
h_{1}:=  K_{1}^{-1}.
\]
By taking a complete orthonormal basis of 
$H^{0}(X,{\cal O}_{X}(2(K_{X}+L) + A)\otimes {\cal I}(h_{L}h))$ with 
respect to  the inner product
\[
(\sigma ,\tau ):= (\sqrt{-1})^{n^{2}}\int_{X}\sigma\wedge\bar{\tau}\cdot h_{1},
\]
we define $K_{2}$ and the singular hermitian metric 
\[
h_{2} := K_{2}^{-1}.
\]
Suppose that we have already constructed $\{ h_{1},\cdots ,h_{m-1}\}$.
We set 
\[
V_{m}:= H^{0}(X,{\cal O}_{X}(m(K_{X}+L) + A)\otimes {\cal I}(h_{L}\cdot h^{m-1}))
\]
By taking a complete orthonormal basis of 
$V_{m}$ with 
respect to  the inner product
\[
(\sigma ,\tau ):= (\sqrt{-1})^{n^{2}}\int_{X}\sigma\wedge\bar{\tau}\cdot h_{L}\cdot h_{m-1},
\]
we define $K_{m}$ and the singular hermitian metric 
\[
h_{m} := K_{m}^{-1}
\]
in the same manner. 
We note that for every $x\in X$ and $m$,
\[
h_{m}^{-1}(x)= K_{m} = \sup \{ \mid\sigma\mid^{2}(x) ;  
\sigma\in V_{m}, 
\int_{X}h_{m-1}\cdot \mid\sigma\mid^{2} = 1\}
\]
holds by definition (cf. \cite[p.46, Proposition 1.4.16]{kr}).

We set  
\[
\nu : = \limsup_{m\rightarrow\infty}\frac{\log \dim H^{0}(X,{\cal O}_{X}(m(K_{X}+L) + A)\otimes{\cal I}(h^{m}))}{\log m},
\]
i.e., $\nu$ is the asymptotic  Kodaira dimension of $(K_{X}+L,dV^{-1}\cdot h)$ (cf. Definition \ref{numerical Kodaira}) and is an integer between $0$ and $n = \dim X$ by Theorem \ref{nak}. 
The following theorem is the main result in this section. 
\begin{theorem}\label{dynamical}(cf. \cite{tu6})
Let $X$, $(L,h_{L})$, $\{K_{m}\}_{m=1}^{\infty}$ and $\{ h_{m}\}_{m=1}^{\infty}$be as above. \\  
\noindent Then 
\[
K_{\infty}:=  \mbox{\em the upper envelope of}\,\,\,\limsup_{m\rightarrow\infty}\sqrt[m]{(m!)^{-\nu}K_{m}}
\]
exists and 
\[
h_{\infty}:= 1/K_{\infty}
\]
is an AZD of $(K_{X} + L,dV^{-1}\cdot h_{L})$, where $dV$ is an arbitrary
$C^{\infty}$ volume form on $X$\footnote{In such a case, 
it may be appropriate to write $K_{X} + (L,h_{L})$ instead of 
$(K_{X} +L,dV^{-1}\cdot h_{L})$}. $\square$ 
\end{theorem}
\subsection{Proof of Theorem \ref{dynamical}}\label{DyAZD} 
By the assumption there exists a positive constant $C$ such that 
\[
h^{0}(X,{\cal O}_{X}(m(K_{X} +L)+ A)\otimes{\cal I}(h^{m})) \leqq C\cdot m^{\nu}
\]
holds for every $m \geqq 1$.

Let us fix a K\"{a}hler form $\omega$ on $X$. 
Let $dV$ be the volume form on $X$ with respect to $\omega$ and 
let $h_{L,0}$ be a $C^{\infty}$ hermitian metric on $L$. 
The following estimate is an easy consequence of the submeanvalue property 
of plurisubharmonic functions.

\begin{lemma}\label{upper}
There exists a positive constant $\tilde{C}$ such that 
for every $m\geqq 1$,
\[
K_{m}\leqq \tilde{C}^{m}\cdot(m!)^{\nu}\cdot (dV)^{m}\cdot 
h_{A}^{-1}\cdot h_{L,0}^{-m}
\]
holds.

\end{lemma}
{\bf Proof}. 
Let $p\in X$ be an arbitrary point. 
Let $(U,z_{1},\ldots ,z_{n})$ be a local cooordinate around $p$ 
such that 
\begin{enumerate}
\item $z_{1}(p) = \cdots = z_{n}(p) = 0$,
\item $U$ is biholomorphic to the open unit polydisk in 
$\mbox{\bf C}^{n}$ with center $O\in \mbox{\bf C}^{n}$ 
by the coordinate,
\item $z_{1},\ldots ,z_{n}$ are holomorphic on a neighbourhood 
of the closure of $U$,
\item there exists a holomorphic frame $\mbox{\bf e}$ of $A$ on the closure of $U$.
\end{enumerate}
Taking $U$ sufficiently small we may assume that there exist  holomorphic 
frames $\mbox{\bf e}_{L}$ of $L$  and  $\mbox{\bf e}_{A}$ of $A$ on $U$
respectively. 
We set 
\[
\Omega := (-1)^{\frac{n(n-1)}{2}}(\sqrt{-1})^{n}dz_{1}\wedge\cdots\wedge dz_{n}\wedge d\bar{z}_{1}\wedge
\cdots\wedge d\bar{z}_{n}.
\]
For every $m\geqq 0$, we set 
\[
B_{m} := \sup_{x\in U}\,\,\frac{K_{m}}{\Omega^{m}\cdot \mid\mbox{\bf e}_{L}\mid^{2m}\cdot \mid\mbox{\bf e}_{A}\mid^{2}}(x).
\]
We note that for any $x\in X$
\[
K_{m}(x) = \sup\{ \mid\phi\mid^{2}(x) ; 
\phi\in V_{m}, 
\int_{X}h_{m-1}\mid\phi\mid^{2} = 1\}
\]
holds.
Let $\phi_{0}$ be the element of $V_{m}$ 
such that 
\[
\int_{X}h_{m-1}\mid\phi_{0}\mid^{2} = 1.
\]
Then there exists a holomorphic function $f$ on $U$ such that 
\[
\phi_{0}\mid_{U} = f\cdot (dz_{1}\wedge\cdots dz_{n})^{m}\cdot \mbox{\bf e}_{L}^{m}\cdot\mbox{\bf e}_{A}
\]
holds.
Then 
\[
\int_{U}\mid\phi_{0}\mid^{2}\cdot h_{A}\cdot h_{L}^{m}\cdot \Omega^{-(m-1)}
= \int_{U}\mid f\mid^{2}\cdot h_{L}(\mbox{\bf e}_{L},\mbox{\bf e}_{L})\cdot h_{A}(\mbox{\bf e}_{A},\mbox{\bf e}_{A})\cdot\Omega
\]
holds.  
On the other hand by the definition of $B_{m-1}$ we see that
\[
\int_{U}h_{A}\mid\phi_{0}\mid^{2}\Omega^{-(m-1)} 
\leqq B_{m-1}\int_{U}h_{m-1}\mid\phi_{0}\mid^{2}
\leqq B_{m-1}
\]
hold.
Combining above inequalities we have that
\[
\int_{U}\mid f\mid^{2}h_{A}(\mbox{\bf e},\mbox{\bf e})\,\Omega
\leqq B_{m-1}
\]
holds.
Let $0 < \delta << 1$ be a sufficiently small number. 
Let $U_{\delta}$ be the inverse image of 
\[
\{ (y_{1},\ldots ,y_{n})\in \mbox{\bf C}^{n}; 
\mid y_{i}\mid < 1 - \delta \}
\]
by the coordinate $(z_{1},\ldots ,z_{n})$.

Then by the subharmonicity of $\mid f\mid^{2}$, there exists a positive constant $C_{\delta}$ independent of $m$ such that 
\[
\mid f(x)\mid^{2} \leqq C_{\delta}\cdot B_{m-1}
\]
holds for every $x\in U_{\delta}$.
Then we have that 
\[
K_{m}(x) \leqq C_{\delta}\cdot B_{m-1}\cdot\mid\mbox{\bf e}_{A}\mid^{2}\cdot \mid\mbox{\bf e}_{L}\mid^{m}\cdot \Omega^{m}(x)
\]
holds for every $x\in U_{\delta}$. 
Summing up the estimates for the orthonormal basis, moving $p$, by the compactness of $X$ we see that 
there exists a positive constant $\tilde{C}$ such that 
\[
K_{m} \leqq \tilde{C}^{m}\cdot (m!)^{\nu}\cdot (dV)^{m}\cdot h_{A}^{-1}\cdot h_{L,0}^{-m} 
\]
holds on $X$.
This completes the proof of Lemma \ref{upper}. \vspace{5mm} $\square$\\

\noindent Let $V$ be a  $\nu$ dimensional nonsingular subvariety such that  
\begin{enumerate}
\item $V$ is not contained in the pluripolar set of $h$. 
\item $(K_{X}+L,h)\mid V$ is big, i.e.,
\[
\limsup_{m\rightarrow \infty} m^{-\nu}\cdot h^{0}(V,{\cal O}_{V}(m(K_{X}+L))\otimes {\cal I}(h^{m}\mid_{V})) > 0
\] 
holds. 
\end{enumerate}
Since $(K_{X}+ L\mid_{V},h\mid_{V})$ is big,
by Theorem \ref{kodaira}, there exists a singular hermitian metric 
$h_{V}$ on $K_{X}+L\mid_{V}$ such that 
\begin{enumerate}
\item $\Theta_{h_{V}}$ is strictly positive 
everywhere on $V$,
\item  
$h\mid_{V} \leqq h_{V}$ 
holds on $V$. 
\end{enumerate}

\noindent Suppose that for some $m\geqq 2$ 
\[
K_{m-1}(x) \geqq C_{m-1}\cdot h_{A}^{-1}\cdot h_{V}^{-m}
\]
holds on $x \in V$. 
Let us estimate $K_{m}(x)$ from below  at the point $x$.
\begin{lemma}\label{lower}
Let $\varepsilon$ be a positive number. 
There exists a positive constant $C(\varepsilon )$ depending 
only on $\varepsilon$ and $V$ such that 
\[
K_{m}(x) \geqq C(\varepsilon )^{m}\cdot (m!)^{\nu}\cdot h_{A}^{-1}\cdot 
(h_{V}^{\varepsilon}\cdot h^{(1-\varepsilon )})^{-m}  
\]
holds at $x$.  
\end{lemma} 
{\bf Proof}.
We shall give a proof only for the case $\varepsilon = 1$. 
The general case follows from the proof below in the same manner.  
We note that 
\[
K_{m}(x) = \sup \{ \mid\sigma\mid^{2}(x) ;  
\sigma\in \Gamma (X,{\cal O}_{X}(A+m(K_{X}+L))), 
\int_{X}h_{m-1}\cdot \mid\sigma\mid^{2} = 1\}.
\]
holds. 
Then by the asymptotic expansion of Bergman kernels (Theorem \ref{MA}), we can extend  any element $\sigma_{x}$ of the fiber $(m(K_{X}+L)+A)_{x}$
of the line bundle of $m(K_{X}+L) + A$ to an element 
\[
\sigma_{V}\in H^{0}(U\cap V,{\cal O}_{V}(m(K_{X}+L) +A)\otimes {\cal I}(h_{L}\cdot h_{V}^{m-1})) 
\]   
with 
\[
\int_{U\cap V}h_{V}^{m}\cdot h_{A}(\sigma_{V},\sigma_{V})\mid dz_{1}\wedge\cdots\wedge dz_{\nu}\mid^{2}
\leqq C\cdot m^{-\nu}\cdot h_{V}^{m}\cdot h_{A}(\sigma_{x},\sigma_{x}) 
\]
where $C$ is a positive constant depending only on $x$ and $V$.
Since $h\mid V \leqq h_{V}$ holds on $V$, this implies that 
\[
\int_{U\cap V}h^{m}\cdot h_{A}(\sigma_{V},\sigma_{V})\mid dz_{1}\wedge\cdots \wedge dz_{\nu}\mid^{2} 
\leqq C\cdot m^{-\nu}\cdot h_{V}^{m}\cdot h_{A}(\sigma_{x},\sigma_{x}) 
\] 
holds. 
Then we can extend $\sigma_{V}$ to 
\[
\sigma_{U} \in H^{0}(U,{\cal O}_{X}(m(K_{X}+L) + A)\otimes {\cal I}(h^{m-1}))
\]
by the $L^{2}$-extension theorem (Theorems \cite{o-t} and \ref{extension}) so that 
\[
(\sqrt{-1})^{n^{2}}\int_{U}\sigma_{U}\wedge\bar{\sigma}_{U}\cdot h^{m-1}\cdot h_{A}
\leqq C_{U}\cdot C\cdot m^{-\nu}\cdot h_{V}^{m}\cdot h_{A}(\sigma_{x},\sigma_{x}) 
\]
where $C_{U}$ is a positive constant depending only on $U$. 
Let $\rho$ be a $C^{\infty}$ function such that 
\begin{enumerate}
\item $\mbox{Supp}\, \rho \subset\subset U$,
\item $0\leqq \rho \leqq 1$ on $X$,
\item $\rho\equiv 1$ on some  neighbourhood $W$ of $x$.
\end{enumerate}
Taking $A$ sufficiently ample, we may assume that 
\[
\sqrt{-1}\partial\bar{\partial}(n\rho \log \sum_{i=1}^{n}\mid z_{i}\mid^{2}) 
+ \Theta_{h_{A}}
\]
is strictly positive on $X$.  
Then we may solve the $\bar{\partial}$-equation
\[
\bar{\partial}u = \bar{\partial}(\rho\cdot \sigma_{U})
\]
with
\[
\int_{X}\exp (-(n+1)\rho)\cdot\log \sum_{i=1}^{n}\mid z_{i}\mid^{2})\mid u\mid^{2}h_{A}\cdot h^{m-1} 
\]
\[
\hspace{30mm}\leqq 
C^{\prime}_{U}\cdot \int_{U}\exp (-(n+1)\rho)\cdot\log \sum_{i=1}^{n}\mid z_{i}\mid^{2})\cdot\mid\bar{\partial}(\rho\cdot \sigma_{U})\mid^{2} dV 
\]
holds, where $\mid\bar{\partial}(\rho\cdot\sigma_{U})\mid^{2}$
denotes the norm with respect to $h_{A}\cdot h_{L}^{m-1}$ and $\omega$ 
and $C^{\prime}_{U}$ is a positive constant depending only on the supremum 
of the norm of $\bar{\partial}\rho$ with respect to $\omega$.
This implies that 
\[
u(x) = 0
\]
and  there exists a positive constant $C$ independent of $(L,h_{L})$ and 
$\sigma_{x}$ such that 
\[
\int_{X}h_{A}\cdot h_{L}\cdot\mid u\mid^{2}
\leqq C\cdot (dV^{-1}\cdot h_{A}\cdot h_{L})(\sigma_{x},\sigma_{x})
\]
holds.
Then 
\[
\sigma := \rho\cdot\sigma_{U} - u
\in H^{0}(X,{\cal O}_{X}(m(K_{X}+L)+A)\otimes{\cal I}(h^{m-1}))
\]
is an extension of $\sigma_{x}$ 
such that 
\[
(\sqrt{-1})^{n^{2}}\int_{X}\sigma\wedge\bar{\sigma}\cdot h^{m-1}\cdot h_{A}
\leqq   C^{\prime}\cdot C_{U}\cdot C\cdot m^{-\nu}\cdot (h_{V}^{m}\cdot h_{A})(\sigma_{x},\sigma_{x}) 
\]
where $C^{\prime}$ is a positive constant depending $V$ and $U$. 
Hence by the extremal property of Bergman kernels,
we obtain that 
\[
K_{m}(x) \geqq   (C^{\prime}\cdot C_{U}\cdot C)^{-1}\cdot m^{\nu}\cdot C_{m-1}\cdot h_{A}^{-1}\cdot 
(h_{V})^{-m} 
\]
holds.   
We note that for every $0 < \varepsilon < 1$,
 $h_{V}^{\varepsilon}\cdot h^{1-\varepsilon}\mid_{V}$ is a singular hermitian metric with strictly positive curvature on $V$. 
Then replacing $h_{V}$ by   $h_{V}^{\varepsilon}\cdot h^{1-\varepsilon}\mid_{V}$we obtain that there exists a positive constant $C(\varepsilon )$ depending 
on $V$ and $\varepsilon$ such that 
\[
K_{m}(x) \geqq  C(\varepsilon )\cdot m^{\nu}\cdot C_{m-1}\cdot h_{A}^{-1}\cdot 
(h_{V}^{\varepsilon}\cdot h^{(1-\varepsilon )})^{-m} 
\]
holds. 
Hence summing up the estimates for $m$, we get 
\[
K_{m}(x) \geqq C(\varepsilon )^{m}\cdot (m!)^{\nu}\cdot h_{A}^{-1}\cdot 
(h_{V}^{\varepsilon}\cdot h^{(1-\varepsilon )})^{-m}  
\]
at $x$.   
The estimate is valid on a neighbourhood of $x$ in $V$ and moving $x$ on $V$,
we may assume that the above estimate is valid on $V$. 

Hence on $V$, there exists a positive constant $C_{V}$ such that 
\[
K_{m}(y) \geqq C_{V}(\varepsilon)^{m}\cdot (m!)^{\nu}\cdot h_{A}^{-1}\cdot 
(h_{V}^{\varepsilon}\cdot h^{(1-\varepsilon )})^{-m}
\]  
holds for every $y\in V$.

Combining Lemmas \ref{upper} and \ref{lower} we have that  the limit :
\[
\limsup_{m\rightarrow\infty}\sqrt[m]{(m!)^{-\nu}K_{m}} 
\]
exists and nonzero on $V$.  Moving $V$, the limit exists on the whole $X$.
We set 
\[
h_{\infty} := \mbox{the lower envelope of}\,\,\,(\limsup_{m\rightarrow\infty}\sqrt[m]{(m!)^{-\nu}K_{m}})^{-1} 
\]

Let us prove $h_{\infty}$ is an AZD of $(K_{X} + L,dV^{-1}\cdot h_{L})$.
Let $x\in X$ be an arbitrary point.  Then there exists a family 
$\{ V_{t}\} (t \in \Delta^{n-\nu})$ of smooth projective subvariety of dimension $\nu$ in $X$ and a local coordinate $(U,z_{1},\cdots ,z_{n})$ such that 
\begin{enumerate}
\item $(z_{1}(x), \cdots ,z_{n}(x)) = O$, 
\item $U$ is biholomorphic to $\Delta^{n}$ via $(z_{1},\cdots ,z_{n})$, 
\item $V_{t}\cap U = \{ p\in U\mid (z_{\nu+1}(p), \cdots ,z_{n}(p)) = t\}$.
\item $(K_{X} + L,h)\mid_{V_{t}}$ is big for every $t\in \Delta^{n-\nu}$. 
\end{enumerate} 
Let 
\[
\phi : V \longrightarrow \Delta^{n-\nu}
\]
be the above family.  
By the proof of Theorem \ref{kodaira}, we may take a singular hermitian metric
$h_{V}$ on $\phi^{*}(K_{X} + L)$ such that
\begin{enumerate}
\item $h_{V}\mid V_{t}$ is a singular hermitian metric with strictly positive
curvature. 
\item There exists a positive constant $C$ such htat 
\[
h\mid_{V_{t}}\leqq C\cdot h_{V}\mid_{V_{t}}
\]
holds on $V_{t}$ for every $t\in \Delta^{n-\nu}$.
\end{enumerate}
Since  for every $y\in U$, there exists a positive constant $C(\varepsilon)$ such that  
\[
K_{m}(y) \geqq C(\varepsilon )^{m}\cdot (m!)^{\nu}\cdot h_{A}^{-1}\cdot 
(h_{V}^{\varepsilon}\cdot h^{(1-\varepsilon )})^{-m}  
\]
holds for every $m \geqq 1$ by Lemma \ref{lower}. 
Then since $\limsup_{m\rightarrow\infty}\sqrt[m]{(m!)^{-\nu}K_{m}}$ exist, letting $m$ tend to infinity, see that 
\[
h_{\infty}\mid_{V_{t}} \leqq C(\varepsilon)^{-1}\!\!\cdot h_{V}^{\varepsilon}\cdot h^{(1-\varepsilon )}
\]
holds. 
Let $\sigma \in \Gamma (U,{\cal O}_{X}(m(K_{X}+L))\otimes {\cal I}_{\infty}(h^{m}))$
be any element. 
Then by the above iequality, we see that 
\[
\int_{U}\mid\sigma\mid^{p}\cdot (h_{\infty})^{p}\,\,dV 
\leqq C(\varepsilon )^{-p}\!\!\cdot \int_{U}\mid\sigma\mid^{p}\cdot (h_{V}^{\varepsilon}h^{(1-\varepsilon )})^{p}\,\, dV
\]
holds. Hence letting $\varepsilon$ tend to $0$ and shrinking $U$, if necessary,  we see that 
\[
\sigma \in L^{p}_{loc}(m(K_{X}+L),h_{\infty}^{p})
\]
holds for every $p\geqq 1$.  
This implies that 
\[
\sigma \in \Gamma (U,{\cal O}_{X}(m(K_{X}+L))\otimes \bar{\cal I}_{\infty}(h_{\infty}^{m}))
\]
holds. 
Hence we see that $h_{\infty}$ an AZD of $(K_{X}+L,dV^{-1}\cdot h_{L})$
(cf. Definition \ref{singAZD}).  
This completes the proof of Theorem \ref{dynamical}. $\square$  

\section{Proof of Theorem  \ref{subad2}}

In this section, we shall prove Theorem \ref{subad2} 
together with the following theorem which is a little bit weaker 
than Theorme \ref{main}.
 
\begin{theorem}\label{subad1}(\cite[Theorem 5.1]{tu5})
Let $X$,$S$,$\Psi_{S}$ be as in Section 1.2.  
Suppose that $S$ is smooth. 
Let $d$ be a positive integer such that $d > \alpha m_{0}$. 
We assume  that $(K_{X}+L,e^{-\varphi}\cdot dV^{-1}\cdot h_{L}\mid_{S})$
 is weakly pseudoeffective (cf. Definition \ref{wpe}) and  let  $h_{S}$ be an AZD of $(K_{X}+L,e^{-\varphi}\cdot dV^{-1}\cdot h_{L}\mid_{S})$.
Suppose that $h_{S}$ is normal (cf. Definition \ref{normal}) or
for every ample line bundle $A$ on $X$ 
\[
\dim H^{0}(S,{\cal O}_{S}(A+ m(K_{X}+L)\mid_{S})\otimes{\cal I}(h_{S}^{m}))
= O(m^{\nu})
\]
holds, where $\nu$ denotes the numerical Kodaira dimension 
$\nu_{num}(K_{X}+L\mid_{S},h_{S})$ of $(K_{X}+L\mid_{S},h_{S})$ (cf. Definition \ref{numerical}).
  
Then every element of 
\[
H^{0}(S,{\cal O}_{S}(m(K_{X}+dL))\otimes{\cal I}(e^{-\varphi}\cdot h_{L}^{d}\cdot h_{S}^{m-1}))
\]
extends to an element of 
\[
H^{0}(X,{\cal O}_{X}(m(K_{X}+dL))\otimes {\cal I}(h_{L}^{d}\cdot h_{0}^{m-1})), 
\]
where $h_{0}$ is an AZD of $K_{X} + dL$ with minimal singularities. 
In particular every element of 
\[
H^{0}(S,{\cal O}_{S}(m(K_{X}+dL))\otimes{\cal I}(e^{-\varphi}\cdot h_{L}^{d})\cdot {\cal I}_{\infty}(h^{m-1}_{L}\mid_{S}\cdot e^{-(m-1)\varphi}))
\] 
extends to an element of 
\[
H^{0}(X,{\cal O}_{X}(m(K_{X}+dL))\otimes {\cal I}(h_{0}^{m})). 
\]
$\square$ \end{theorem}
\begin{remark}
The normality condition on $h_{S}$ is more restrictive than the condition
\[
\nu_{num}(K_{X}+L\mid_{S},h_{S}) = \nu_{asym}(K_{X}+L\mid_{S},h_{S})
\]
(see Theorem \ref{nak}).  But practically, the normality of $h_{S}$ is
much easier to verify than the equality of $\nu_{num}$ and $\nu_{asymp}$. 
  For example  if $h_{L}\!\!\mid\!S$ is normal and $K_{X}\!\!\mid\! S$ is $\mathbb{Q}$-effective or admits a normal singular hermitian metric with semipositive curvature current, then $h_{S}$ is normal. $\square$ \vspace{5mm}  
\end{remark}
\subsection{Setup}
Let $X$ be a smooth projective variety 
and let $(L,h_{L})$ be a singular hermitian line bundle on $X$ such that 
$\Theta_{h_{L}}\geqq 0$ on $X$.  
As we mensioned in Section 1.2, we assume that $h_{L}$ is lowersemicontinuous. 
This is a technical assumption so that a local potential 
of the curvature current of $h_{L}$ is plurisubharmonic. 

Let $m_{0}$ be a positive integer and 
let $\sigma \in \Gamma (X,{\cal O}_{X}(m_{0}L)\otimes {\cal I}(h^{m_{0}}_{L}))$ be a 
global section. 
Let $\alpha$ be a positive rational number $\leqq 1$ and let $S$ be 
an irreducible subvariety of $X$ 
such that  $(X, \alpha (\sigma ))$ is LC(log canonical) but not KLT(Kawamata log terminal)
on the generic point of $S$ and $(X,(\alpha -\epsilon )(\sigma ))$ is KLT on the generic point of $S$ 
for every $0 < \epsilon << 1$. 
We set 
\[
\Psi_{S} = \alpha \log h_{L}(\sigma ,\sigma ).
\]
Suppose that $S$ is smooth for simplicity 
(if $S$ is not smooth, we just need to take an embedded 
resolution to apply Theorem \ref{subad1}). 
We shall assume that $S$ is not contained in the 
singular locus of $h_{L}$, where the singular locus of $h_{L}$ means the 
set of points where $h_{L}$ is $+\infty$. 
Let $dV$ be a $C^{\infty}$ volume form on $X$. 

Then  we may define a (possibly singular) measure 
$dV[\Psi_{S}]$ on $S$ as in the introdcution. 
Let $dV_{S}$ be a $C^{\infty}$ volume form on $S$ and 
let $\varphi$ be the function on $S$ defined by
\[
\varphi := \log \frac{dV_{S}}{dV[\Psi_{S}]}
\]
($dV[\Psi_{S} ]$ may be singular on a subvariety of $S$, also 
it may be totally singular on $S$). 

\subsection{Dynamical construction of singular hermitian metrics
with successive extensions}\label{Dycon}
Let us start the proof of Theorem \ref{subad1}.
Let $d$ be a positive integer such that $d > \alpha m_{0}$. 
Replacing $L$ by $dL$, we may and do assume that $d = 1$ from the 
beginning. 
By the assumption, there exists an AZD   $h_{S}$  of \\ $(K_{X} + L\mid_{S}, (dV^{-1}\cdot h_{L})\mid_{S}\cdot e^{-\varphi})$.   
We shall define  sequences of the hermitian metrics 
$\{ h_{m}\}$ on $A + m(K_{X}+L)\mid S$ and $\{ \tilde{h}_{m}\} (m\geqq 1)$ 
on $A +m(K_{X}+L)$ inductively as follows.
Let $h_{A}$ be a $C^{\infty}$ hermitian metric on $A$.
Let 
\[
\{ \sigma^{(1)}_{0},\cdots ,\sigma^{(1)}_{N_{1}}\}
\]
be an orthonormal basis of $H^{0}(S,{\cal O}_{S}(A+(K_{X}+L))\otimes {\cal I}(h_{S}))$ with respect to the inner product :
\[
(\sigma ,\tau ) = \int_{S}\sigma\cdot \bar{\tau}\cdot (h_{A}\cdot h_{L}\cdot dV^{-1})dV[\Psi_{S}] (\sigma ,\tau \in  H^{0}(S,{\cal O}_{S}(A+(K_{X}+L))\otimes {\cal I}(h_{S}))). 
\]
We set 
\[
K_{1} = \sum_{i=0}^{N_{1}}\mid\sigma^{(1)}_{i}\mid^{2}
\]
and set 
\[
h_{1} = 1/K_{1}. 
\]
Then $h_{1}$ is a singular hermitian metric on $A+(K_{X}+L)\mid S$. 
By the choice of $A$, we see that ${\cal O}_{S}(A+(K_{X}+L))\otimes {\cal I}(h_{S})$ is globally genegerated on $S$. 
Hence we see that 
\[
h_{1} \leqq O(h_{A}\cdot h_{S}) 
\]
holds, where this means that $h_{1}\cdot (h_{A}\cdot h_{S})^{-1}$ is 
bounded from above on $S$. 
Then by the $L^{2}$-extension theorem (Theorem \ref{extension}), each $\sigma_{i}^{(1)}$ extends to a section
\[
\tilde{\sigma}^{(1)}_{i} 
\in H^{0}(X,{\cal O}_{X}(A+(K_{X}+L)))
\]
such that 
\[
\parallel\sigma_{i}^{(1)}\parallel 
= (\int_{X}\mid\sigma_{i}^{(1)}\mid^{2}\cdot (h_{A}\cdot h_{L}\cdot dV^{-1})
\cdot dV)^{\frac{1}{2}}
\]
satisfies the inequality 
\[
\parallel\sigma_{i}^{(1)}\parallel \leqq C, 
\]
where $C$ is the positive constant as in Theorem \ref{extension}.
And we set 
\[
\tilde{K}_{1} = \sum_{i=0}^{N_{1}} \mid\sigma_{i}^{(1)}\mid^{2}
\]
and
\[
\tilde{h}_{1} := 1/\tilde{K}_{1}.
\]
We note that $\tilde{K}_{1}$ depends on the choice of the 
orthonormal basis.

Suppose that we have already constructed 
$\{ h_{i}\}_{i=1}^{m}$ and $\{\tilde{h}_{i}\}_{i=1}^{m}$ and 
\[
h_{i} = O(h_{A}\cdot h_{S}^{i})
\] 
holds for every $i= 1,\cdots ,m$. 

Let 
\[
\{ \sigma^{(m+1)}_{0},\cdots ,\sigma^{(m+1)}_{N_{m+1}}\}
\]
be an orthonormal basis of $H^{0}(S,{\cal O}_{S}(A+(m+1)(K_{X}+L))\otimes {\cal I}(h^{m+1}_{S}))$ with respect to the inner product :
\[
(\sigma ,\tau ) = \int_{S}\sigma\cdot \bar{\tau}\cdot (h_{m}\cdot h_{L}\cdot dV^{-1})dV[\Psi_{S}].
\]
Here we note that by the assumption $h_{m} = O(h_{A}\cdot h_{S}^{m})$, 
this inner product is well defined. 
By the choice of $A$, we see that 
${\cal O}_{S}(A + (m+1)(K_{X}+L))\otimes {\cal I}(h_{A}h_{S}^{m+1})$ is 
globally generated.  Hence 
\[
h_{m+1} = O(h_{A}\cdot h_{S}^{m+1}) 
\]
holds. 

We define
\[
K_{m+1} := \sum_{i=0}^{N_{m+1}}\mid\sigma^{(m+1)}_{i}\mid^{2} 
\]
and 
\[
h_{m+1} := 1/K_{m+1}.
\]
Then $h_{m+1}$ is a singular hermitian metric on $A + (m+1)(K_{X}+L)$
with semipositive curvature current. 
 
Then by the $L^{2}$-extension theorem (Theorem \ref{extension}), each $\sigma_{i}^{(m+1)}$ extends to a section
\[
\tilde{\sigma}^{(m+1)}_{i} 
\in H^{0}(X,{\cal O}_{X}(A+(K_{X}+L)))
\]
such that 
\[
\parallel\tilde{\sigma}_{i}^{(m+1)}\parallel 
= (\int_{X}\mid\tilde{\sigma}_{i}^{(m+1)}\mid^{2}\cdot (h_{L}\cdot dV^{-1}\cdot \tilde{h}_{m})
\cdot dV)^{\frac{1}{2}}
\]
satisfies the inequality 
\[
\parallel\tilde{\sigma}_{i}^{(m+1)}\parallel \leqq C, 
\]
where $C$ is the positive constant independent of $m$ as in Theorem \ref{extension}. 
And we set 
\[
\tilde{K}_{m+1} := \sum_{i=0}^{N_{m+1}} \mid\sigma_{i}^{(m+1)}\mid^{2}.
\]
and 
\[
\tilde{h}_{m+1}:= 1/\tilde{K}_{m+1}. 
\]
\noindent In this way we construct the sequences $\{ K_{m}\}_{m=1}^{\infty}$,
$\{\tilde{K}_{m}\}_{m=1}^{\infty}$, $\{ h_{m}\}_{m=1}^{\infty}$ 
and $\{\tilde{h}_{m}\}_{m=1}^{\infty}$.  Next we shall discuss the 
(normalized) convergence of these sequences. 
Let $\nu$ be the asymptotic  Kodaira  dimension (cf. Definition \ref{numerical Kodaira}) of $(K_{X}+L\mid_{S},h_{S})$, i.e.,
\[
\nu := \limsup_{m\rightarrow\infty}\frac{\log \dim H^{0}(S,{\cal O}_{S}(A + m(K_{X}+L))\otimes {\cal I}(h_{S}^{m}))}{\log m}.
\]
$\nu$ is a nonnegative integer between $0$ and $\dim S$ by Theorem \ref{nak}.  

\begin{lemma}\label{convergence}
\[
h_{\infty}:= \mbox{\em the lower envelope of }\,\,\,\liminf_{m\rightarrow\infty}\sqrt[m]{(m!)^{\nu}\cdot h_{m}}
\]
exists as a singular hermitian metric on $K_{X}+ L\mid_{S}$ and 
is an AZD of  \\ $(K_{X} + L\mid_{S},dV^{-1}\cdot h_{L}\cdot e^{-\varphi})$ on $S$. And 
\[
\tilde{h}_{\infty} : = \mbox{\em the lower envelope of}\,\,\, \liminf_{m\rightarrow\infty}\sqrt[m]{(m!)^{\nu}\cdot \tilde{h}_{m}}
\]
exists as a singular hermitian metric on $K_{X} + L$ with semipositive curvature current and
\[
h_{\infty} \geqq \tilde{h}_{\infty}\!\mid_{S}
\]
holds on $S$. 
$\square$
\end{lemma}
{\bf Proof of Lemma \ref{convergence}.}
To prove Lemma \ref{convergence}, we shall estimate $K_{m}$  from above and below and $\tilde{K}_{m}$
from above. 

The estimate for $\{ K_{m}\}$ is identical as the one in the proof of Theorem \ref{dynamical} in the last section.   
Hence we obtain that 
\[
h_{\infty}:= \mbox{the lower envelope of}\,\,\, \liminf_{m\rightarrow\infty}\sqrt[m]{(m!)^{\nu}\cdot h_{m}}
\]
exists and is an AZD of $(K_{X} + L\mid_{S},dV^{-1}\cdot h_{L}\cdot e^{-\varphi})$ on $S$.

Let us fix $h_{L,0}$ be a $C^{\infty}$ hermitian metric on $L$. 
By the same proof as  that of Lemma \ref{upper}, we obtain that
there exists a positive constant $C_{+}$ such that  
\[
\tilde{K}_{m} \leqq  (C_{+})^{m}\cdot (m!)^{\nu}\cdot h_{A}^{-1}\cdot (dV)^{m}\cdot h_{L,0}^{-m}
\]
holds on $X$ for every $m\geqq 0$. 
Hence by this estimate,  we see that 
\[
\tilde{K}_{\infty}:= \mbox{the upper envelope of}\,\,\,\limsup_{m\rightarrow\infty}\sqrt[m]{\tilde{K}_{m}}
\]
exists on $X$ and is an extension of $K_{\infty}$ by the construction.
In particular $\tilde{K}_{\infty}$ is not identically zero on $X$.
Hence 
\[
h_{\infty} := \frac{1}{\tilde{K}_{\infty}}
\]
is a well defined singular hermitian metric on $K_{X} + L$.
And 
\[
h_{\infty} \geqq \tilde{h}_{\infty}\!\mid_{S}
\] 
holds by the construction. 
This completes the proof of Lemma \ref{convergence}.  $\square$. \vspace{5mm} \\
 
\noindent Let us complete the proof of Theorem \ref{subad1}. 
By Lemma \ref{convergence} and Theorem \ref{extension}, we see that 
every element of 
\[
H^{0}(S,{\cal O}_{S}(m(K_{X} + L))\otimes {\cal I}(e^{-\varphi}\cdot h_{L}\mid_{S}\cdot h_{S}^{m-1}))
\]
extends to an element of 
\[
H^{0}(X,{\cal O}_{X}(m(K_{X} + L)\otimes {\cal I}(h_{L}\cdot\tilde{h}_{\infty}^{m-1})).
\]
Let $h_{0}$ be an AZD of $K_{X} + L$ of minimal singularities. 
Since $h_{\infty}$ has semipositive curvature in the sense of current, 
we see that there exists a positive constant $C_{\infty}$ such that  
\[
h_{0} \leqq C_{\infty}\cdot h_{\infty}
\]
holds. 
Hence  we see that 
every element of 
\[
H^{0}(S,{\cal O}_{S}(m(K_{X} + L))\otimes {\cal I}(e^{-\varphi}\cdot h_{L}\mid_{S}\cdot h_{S}^{m-1}))
\]
extends to an element of 
\[
H^{0}(X,{\cal O}_{X}(m(K_{X} + L)\otimes {\cal I}(h_{L}\cdot h_{0}^{m-1})).
\]
This completes the proof of Theorem \ref{subad1}. $\square$  \vspace{5mm} \\

\noindent The proof of Theorem \ref{subad2} is very similar.
Hence we shall indicate the necessary change.  
The method of the proof is an extension of an AZD of 
$(K_{X} + L + \frac{1}{m}E\mid_{S},dV^{-1}\cdot h_{L}\cdot h_{E}^{\frac{1}{m}})$.   Here we need to perform the dynamical construction of 
AZD on the fractional singular hermitian line bundle.  
The necessary change is that we need to tensorize $(E,h_{E})$ at every $m$ 
step instead of tensorize $(L,h_{L})$  at every step as above (see Section 9.1 below, where one can see the concrete construction). 
$\square$
\section{Proof of Theorem \ref{numsubad}}

In this section we shall prove Theorem \ref{numsubad}.  
The method of the proof is parallel to that of \cite{tu5,tu6} except 
the use of Theorems \ref{subad1} and \ref{subad2}. 
\subsection{Positivity result}
In \cite{ka2}, Y. Kawamata proved the following important theorem. 
\begin{theorem}\label{pos}(\cite[p.894,Theorem 2]{ka2})
Let $f : X \longrightarrow B$ be a surjective morphism of smooth projective 
varieties with connected fibers.
Let $P = \sum P_{j}$ and $Q = \sum_{\ell}Q_{\ell}$ be normal crossing divisors on $X$ and $B$ respectively, such that $f^{-1}(Q) \subset P$ and $f$ 
is smooth over $B\backslash Q$.
Let $D = \sum d_{j}P_{j}$ be a $\mathbb{Q}$-divisor on $X$, where $d_{j}$ may be positive, zero or negative, which satisfies the following conditions :
\begin{enumerate}
\item $D = D^{h} + D^{v}$ such that 
$f :\mbox{Supp}(D^{h})\rightarrow B$ is surjective and smooth over $B\backslash Q$, and $f(\mbox{Supp}(D^{v}))\subset Q$.
An irreducible component of $D^{h}$(resp. $D^{v}$) is called horizontal
(resp. vertical).
\item $d_{j} < 1$ for all $j$.
\item The natural homomorphism ${\cal O}_{B}\rightarrow f_{*}{\cal O}_{X}(\lceil -D\rceil )$ is surjective at the generic point of $B$.
\item $K_{X} +  D\sim_{\mathbb{Q}}f^{*}(K_{B} + L)$ for some 
$\mathbb{Q}$-divisor $L$ on $B$.
\end{enumerate} 
Let 
\begin{eqnarray*}
f^{*}Q_{\ell}& =  &\sum_{j}w_{\ell j}P_{j} \\
\bar{d}_{j} & :=  & \frac{d_{j} +w_{\ell j}-1}{w_{\ell j}}\,\,\,\,\mbox{if}\,\,\,\,
f(P_{j}) = Q_{\ell} \\
\delta_{\ell} &: =  & \max \{\bar{d}_{j} ; f(P_{j}) = Q_{\ell}\} \\
\Delta & :=  & \sum_{\ell}\delta_{\ell}Q_{\ell} \\
M & :=  & L - \Delta .
\end{eqnarray*}
Then $X$ is nef.  $\square$
 \end{theorem} 
\begin{remark} In Theorem \ref{pos}, the condition: 
$d_{j} < 1$ is irrelevant for every $D_{j}$ with $f(D_{j}) \subset Q$
by a trivial reason. 
In fact in this case,  if we replace $D$ by $D^{\prime}:= D - \alpha f^{*}Q$ and replace
$L$ by $L^{\prime}:= L -\alpha Q$ for a sufficiently large positive rational number $\alpha$, $D^{\prime}= \sum d_{j}^{\prime}D_{j}$ satisfies the condition :
$d_{j}^{\prime} < 1$ for all $j$. $\square$ 
\end{remark}
 
Here the meaning of the divisor $\Delta$ may be difficult to understand.
So I would like to give an geometric interpretation of $\Delta$.  
Let $X,P,Q,D,B,\Delta$ be as above. Let $dV$ be a 
$C^{\infty}$ volume form on $X$. 
Let $\sigma_{j}$ be a global section of ${\cal O}_{X}(P_{j})$
with divisor $P_{j}$. 
Let $\parallel\sigma_{j}\parallel$ denote the hermitian norm 
of $\sigma_{j}$ with respect to a $C^{\infty}$ hermitian metric
on ${\cal O}_{X}(P_{j})$ respectively. 
Let us consider the singular volume form
\[
\Omega := \frac{dV}{\prod_{j}\parallel\sigma_{j}\parallel^{2d_{j}}}
\]
on $X$.
Then by taking the fiber integral of $\Omega$ with respect to 
$f : X \longrightarrow B$, we obtain a singular volume form 
$\int_{X/B}\Omega$ on $B$, where the fiber integral $\int_{X/B}\Omega$
is defined by the property that for any open set $U$ in $B$, 
\[
\int_{U}(\int_{X/B}\Omega ) = \int_{f^{-1}(U)}\Omega
\]
holds. 
We note that the  condition 2 in Theorem \ref{pos} assures that 
$\int_{X/B}\Omega$ is continuous on a nonempty Zariski open subset 
of $B$.
Also by the condition 4 in Theorem \ref{pos}, computing the differential 
$df$, we see that 
$K_{X}+D$ is numerically $f$-trivial and 
$(\int_{X/B}\Omega)^{-1}$ is  a $C^{0}$-hermitian metric 
on the $\mathbb{Q}$-line bundle $K_{B}+\Delta$.
Thus  the divisor $\Delta$ corresponds exactly to  
singularities (poles and degenerations)  
of the singular volume form $\int_{X/B}\Omega$ on $B$.

\subsection{Proof of Theorem \ref{numsubad}}

Let $\sigma$ be a multivalued holomorphic section of $D$ 
with divisor $D$.  Let $h$ be the supercanononical AZD (\cite{tu9}) of $K_{X}$
or any AZD constructed in the proof of Theorem \ref{AZD}. 
Let $dV$ be a $C^{\infty}$ volume form on $X$ and let 
\[
\Psi := \log h^{\alpha}(\sigma ,\sigma)
\]
Let $dV[\Psi_{S}]$ be the residue volume form on $S$ defined by 
\[
dV[\Psi] = \mbox{Res}_{S}(e^{-\Psi}dV). 
\]
 
Let 
\[
\pi : Y \longrightarrow X
\]
be a log resolution of  $(X,D)$.
Then by the assumption there exist
irreducible components of $\pi^{-1}(S)$ with discrepancy $-1$ 
which dominates $S$.
We divide the proof into the following two cases. 
\begin{enumerate}
\item There exists a unique irreducible component of $\pi^{-1}(S)$ with discrepancy $-1$ which dominates $S$. 
\item There exist several irreducible components of $\pi^{-1}(S)$ with discrepancy $-1$ which dominate $S$ 
\end{enumerate}
In the first case,  the residue volume form $dV[\Psi_{S}]$ is not identically 
$+\infty$ on $S$.  
In the second case the residue volume form $dV[\Psi_{S}]$ is  identically 
$+\infty$ on $S$.  We shall reduce this case to the first case above by 
a minor modification.     

First we shall consider the first case.
Let 
\[
\varpi : F \longrightarrow S
\]
be the restriction of $\pi$ to $F$.  
We shall write  
\[
\pi^{*}(K_{X}+D) = K_{Y} + F + E,
\]
where $\mbox{Supp} (F + E)$ is a divisor with normal crossings and 
$F$ and $E$ have no common divisorial component.
We set $G := E\mid_{F}$. 
Then we have that there exists a $\mathbb{Q}$-divisor $L$ on $S$.  
\begin{equation}\label{adj}
K_{F} + G =  \varpi^{*}(K_{S} + L). 
\end{equation}
Now we shall apply Theorem  \ref{pos} and obtain that 
\[
L - \Delta
\]
is nef, where $\Delta$ is the $\mathbb{Q}$-divisor on $S$ defined as in 
Theorem \ref{pos}. 
We note that $\Delta$ is effective in this case, since $S$ is smooth. 
Let $\sigma_{\Delta}$ be a multivalued holomorphic section of $\Delta$ 
on $S$ with divisor $\Delta$ and let $h_{\Delta}$ be a $C^{\infty}$ hermitian metric on $\Delta$.
Let $dV_{S}$ be a $C^{\infty}$ volume form on $S$.
Then by the definition of $dV[\Psi]$, there exists a positive constant $C$ such that   
\[
dV[\Psi] \leqq  C\cdot\frac{ dV^{-\alpha}\cdot h^{-\alpha}}{h_{\Delta}(\sigma_{\Delta},\sigma_{\Delta})}\cdot dV_{S}
\]
holds on $S$. We set 
\[
\varphi = \log \frac{dV_{S}}{dV[\Psi]}.
\] 
Then  by the above inequality, we have that 
\begin{equation}\label{key ineq}
e^{-\varphi}\cdot h^{\alpha}\leqq C \cdot dV^{-\alpha}\cdot \frac{1}{h_{\Delta}(\sigma_{\Delta},\sigma_{\Delta})}
\end{equation}
holds.   

Let $d$ be a positive integer greater than $\alpha$. 
By (\ref{adj}),(\ref{key ineq}) and the facts that $L- \Delta$ is nef and $K_{S}$ is pseudoeffective, we see  that
$(1+d)K_{X}\mid_{S}$  
admits a singular hermitian metric with semipositive curvature current which dominates \\
$e^{-\varphi}\cdot(dV\mid_{S})^{-1}\cdot h^{d}$. 
Let $h_{A}$ be a $C^{\infty}$ hermitian metric on $A$ with strictly positive curvature.    
Then 
$((1+d)K_{X}+\frac{1}{m}A\mid_{S},e^{-\varphi}\cdot (dV\mid_{S})^{-1}\cdot h^{d}\cdot h_{A}^{\frac{1}{m}})$  
is big and admits an AZD $h_{S,m}$. 
By Theorem \ref{subad2}, we have the following lemma. 
\begin{lemma}\label{ext}
Let $A$ be an ample line bundle on $X$. 
For every positive integer $m$, every element of 
\[
H^{0}(S,{\cal O}_{S}((m+1)(1 + d)K_{X}+A)\otimes {\cal I}(e^{-\varphi}\cdot h^{d}\cdot h_{S,m}^{m}))
\]
extends to an element of 
\[
H^{0}(X,{\cal O}_{X}((m+1)(1 +d)K_{X}+A)\otimes {\cal I}(h^{(1+d)(m+1)})). 
\]
$\square$ 
\end{lemma}

On the other hand by (\ref{key ineq}), we have that
there exists an inclusion  
\[
H^{0}(S,{\cal O}_{S}((m+1)(1 + d)K_{X}+A)\otimes {\cal I}(h^{(m+1)(d -\alpha)}\cdot 
h_{\Delta}(\sigma_{\Delta},\sigma_{\Delta})^{-1}))
\hookrightarrow 
\]
\[
\hspace{30mm}
H^{0}(S,{\cal O}_{S}((m+1)(1 + d)K_{X}+A)\otimes {\cal I}(e^{-\varphi}\cdot h^{d}\cdot h_{S,m}^{m}))
\]
for every $m \geqq 1$. 
Since $L - \Delta$ is nef and 
\[
\pi^{*}(K_{X} + D)\mid_{F} \sim_{\mathbb{Q}} \varpi^{*}(K_{S} + L) 
\]
holds, if we take a suffienctly ample line bundle $B$ on $X$, 
there exists an inclusion   
\[
H^{0}(S,{\cal O}_{S}((m+1)K_{S}+A))
\hookrightarrow  \hspace{50mm}
\]
\[
\hspace{30mm} H^{0}(S,{\cal O}_{S}((m+1)(1 + d)K_{X}+A+ B)\otimes {\cal I}(h^{(m+1)(d -\alpha)}\cdot 
h_{\Delta}(\sigma_{\Delta},\sigma_{\Delta})^{-(m+1)}))
\]
for every $m \geqq 1$. 
Hence by (\ref{key ineq}), we see that there exists a natural inclusion 
\[
H^{0}(S,{\cal O}_{S}((m+1)K_{S}+A))
\hookrightarrow 
H^{0}(S,{\cal O}_{S}((m+1)(1 + d)K_{X}+A+ B)\otimes {\cal I}(e^{-\varphi}\cdot h^{d}\cdot h_{S,m}^{m}))).
\]
By  the above inclusion  
and  Lemma \ref{ext}, we obtain an injection 
\[
\hspace{-50mm} H^{0}(S,{\cal O}_{S}(mK_{S} + A)) 
\hookrightarrow 
\]
\[
\mbox{Image}\{ H^{0}(X,{\cal O}_{X}(m(1+d)K_{X}+A+B))
\rightarrow H^{0}(S,{\cal O}_{S}(m(1+d)K_{X} + A+B))\}
\]  
for every $m \geqq 1$. 

Next we shall consider the second case, i.e., 
there are several divisorial component of  discrepancy $-1$ of $\pi^{*}(K_{X} +D)$ which dominates $S$.  Let $A$ be an ample divisor on $X$ as above.  
We may assume that $\mbox{Supp}\, R$ contains all the component of 
$\pi^{*}D$.  Let us $(X,D)$ by  $(X,(1 + \delta (\varepsilon ))(D - \varepsilon R))$
where $1 + \delta (\varepsilon)$ is the log canonical threshold of 
$D - \varepsilon R$ along $S$. 
Then perturbing the coefficients of $R$, if necessary, we may assume that 
there exists a unique irredicible divisor with discrepancy $-1$ over $S$ 
Let $\ell$ be a sufficiently positive integer such that 
$\ell A^{\prime}$ is a very ample Cartier divisor and let 
$\{\tau_{0},\cdots ,\tau_{N}\}$ a basis of $H^{0}(X,\pi_{*}{\cal O}_{Y}(\ell A^{\prime}))$.
Then we modify $\Psi$ as 
\[
\Psi_{m} := \frac{1}{1+ \delta (\varepsilon /\ell m)}\{ \log h^{\alpha}(\sigma ,\sigma)
+ \frac{\varepsilon}{\ell m}\log (\sum_{i=0}^{N}h_{A}^{\ell}(\tau_{i},\tau_{i}))\},
\]
where the parameter $\varepsilon$ is a positive number less than $1$. 
Then by the choice of $R$, we see that $dV[\Psi_{m}]$ is not identically $+\infty$ on $S$. 
As Lemma \ref{ext}, for every ample line bundle $A$ on $X$, we can extend every element of    
\[
H^{0}(S,{\cal O}_{S}((m+1)(1 + d)K_{X}+A)\otimes {\cal I}(e^{-\varphi_{m}}\cdot h^{d}\cdot h_{S,m}^{m}))
\]
to an element of 
\[
H^{0}(X,{\cal O}_{X}(m(1 +d)K_{X}+A)\otimes {\cal I}(h^{m})) 
\]
for every $m \geqq 1$, where $h_{S,m}$ is an AZD of 

$((1 + d)K_{X} + \frac{1}{m}A,e^{-\varphi_{m}}\cdot (dV\mid_{S})^{-1}\cdot h^{d}\cdot h_{A}^{\frac{1}{m}})$ 
Then replacing $\varphi$ by $\varphi_{m}$ and tracing the proof of the first case, if we take a sufficiently ample line bundle $B$ on $X$  and take $\varepsilon$ sufficiently small, again we obtain an injection
\[
\hspace{-50mm} H^{0}(S,{\cal O}_{S}(mK_{S} + A)) 
\hookrightarrow 
\]
\[
\mbox{Image}\{ H^{0}(X,{\cal O}_{X}(m(1+d)K_{X}+A+B)\otimes {\cal I}(h^{m}))
\rightarrow H^{0}(S,{\cal O}_{S}(m(1+d)K_{X} + A +B))\}
\]  
for every $m\geqq 1$. 
Hence we obtain that 
\[
\nu ((K_{X},h)\mid_{S}) \geqq \nu (K_{S})
\]
holds. 
This completes the proof of Theorem \ref{numsubad}.  $\square$

\section{Proof of Theorems  \ref{main} and \ref{numsubad2}}

Let $M,S,(L,h_{L}),\Psi_{S},dV,\varphi$ be as in Theorem \ref{main}. 
By the assumption \\ $(K_{X} + L\mid_{S},e^{-\varphi}\cdot dV^{-1}\cdot h_{L}\mid_{S})$ is weakly pseudoeffective.   
We set $n:= \dim S$.

\subsection{Dynamical construction with a parameter}
Let $A$ be an ample line bundle on $X$ and let $h_{A}$ be a $C^{\infty}$ hermitian metric on $A$ with strictly positive curvature.
Then we see that $(K_{X} + L + \frac{1}{\ell}A\mid_{S},e^{-\varphi}\cdot dV^{-1}\cdot h_{L}\cdot h_{A}^{\frac{1}{\ell}}\mid_{S})$   is big for every $\ell > 0$. 
The main idea of the proof of Theorem \ref{main} is to extend 
an AZD of  $(K_{X} + L + \frac{1}{\ell}A\mid_{S},e^{-\varphi}\cdot dV^{-1}\cdot h_{L}\cdot h_{A}^{\frac{1}{\ell}}\mid_{S})$ to a singular hermitian metric
of $K_{X} + L + \frac{1}{\ell}A$ of semipositive curvature current with uniform estimates with respect to 
the parameter $\ell$.   And prove the normalized convergence 
as $\ell$ tends to infinity.  

Let $h_{S}$ be an AZD of $(K_{X} + L\mid_{S},e^{-\varphi}\cdot dV^{-1}\cdot h_{L}\mid_{S})$ with minimal singularities. 
Then as in Section 7.2, we construct an AZD on $(K_{X} + L + \frac{1}{\ell}A\mid_{S},e^{-\varphi}\cdot dV^{-1}\cdot h_{L}\mid_{S})$ by using the dynamical system of Bergman kernels as 
\[
K_{\ell ,1} := \sum_{i}\mid\!\sigma_{i}^{(\ell,1)}\!\!\mid^{2},
\]
where $\{\sigma_{0}^{(\ell,1)},\cdots ,\sigma_{N(\ell,1)}^{(\ell,1)}\}$
is an  orthonormal basis of 
\[
H^{0}(S,{\cal O}_{S}(K_{X}+L + A\mid_{S})\otimes {\cal I}(h_{S}))
\]
with respect to the inner product:
\[
(\sigma,\sigma^{\prime}) := \int_{S}\sigma\cdot\bar{\sigma}^{\prime}\cdot dV^{-1}\cdot h_{L}\cdot h_{A}\cdot dV[\Psi_{S}].
\]
And we define 
\[
h_{\ell,1} : = \frac{1}{K_{\ell,1}}. 
\]
Suppose that we have already defined  the singular hermitian metrics  $\{ h_{\ell,1},\cdots ,h_{\ell,m-1}\}$, where $h_{\ell,j}$ ($0\leqq j\leqq m-1$) is 
a singular hermitian metric on \\
$j(K_{X}+L)+(1+\lfloor j/\ell\rfloor )A$ respectively.
Then we define $K_{\ell,m}$ and $h_{\ell,m}$ by 
\[
K_{\ell,m}:= \sum_{i=0}^{N(\ell,m)}
\mid\!\sigma_{i}^{(\ell,m)}\!\!\mid^{2}
\]
where $\{\sigma_{0}^{(\ell,m)},\cdots ,\sigma_{N(\ell,m)}^{(\ell,m)}\}$
is an orthonormal basis of 
\[
H^{0}(S,{\cal O}_{S}(m(K_{X} + L) +(1+ \lfloor\frac{m}{\ell}\rfloor )A\mid_{S})
\otimes {\cal I}(h_{S}^{m}))
\]
with repect to the inner product 
\[
(\sigma,\sigma^{\prime}) := \int_{S}\sigma\cdot\bar{\sigma}^{\prime}\cdot dV^{-1}\cdot h_{L}\cdot h_{A}^{(\lceil\frac{m}{\ell}\rceil -\lceil\frac{m-1}{\ell}\rceil)}\cdot h_{\ell,m-1}\cdot dV[\Psi_{S}].
\]
and we define the singular hermitian metric 
on $m(K_{X} + L) +(1+ \lfloor\frac{m}{\ell}\rfloor )A\mid_{S}$ by 
\[
h_{\ell ,m}:= \frac{1}{K_{\ell,m}}.
\]
In this way we construct the dynamical system of Bergman kernels 
$\{ K_{\ell,m}\}$ and the dynamical system of singular hermitian metrics
$\{ h_{\ell,m}\}$ repsectively with the parameter $\ell$.  
\subsection{Upper estimate}
We set 
\[
d_{\ell ,m} = \dim H^{0}(S,{\cal O}_{S}(m(K_{X}+L) +(1+ \lfloor \frac{m}{\ell}\rfloor ) A)\otimes {\cal I}(h_{S}^{m})). 
\]
We shall fix a $C^{\infty}$ hermitian metric $h_{L,0}$ on $L$.  
In view of the proof of Lemma 6.2, we see that there exists a positive constant $C$ independent of $\ell$ and $m$ such that  
\[
K_{\ell,m+1} \leqq C\cdot d_{\ell,m}\cdot K_{\ell,m}\cdot dV\cdot h_{L,0}^{-1}\cdot h_{A}^{-(\lfloor \frac{m+1}{\ell}\rfloor - \lfloor \frac{m}{\ell} \rfloor )}  \]
holds. 
Then  summing up the estimates, we have that 
\begin{equation}\label{up1}
K_{\ell,m} \leqq C_{0}\cdot C^{m-1}\cdot (\prod_{j=1}^{m-1}d_{\ell,j})\cdot h_{A}^{-(1+\lfloor \frac{m}{\ell}\rfloor)}\cdot dV^{m}\cdot h_{L,0}^{-m}, 
\end{equation}
where $C_{0}$ is a positive constant such that 
\[
K_{\ell,1} \leqq C_{0}\cdot h_{A}^{-1}\cdot dV \cdot h_{L,0}^{-1}
\]
holds on $S$.  
We note that by the definition of $(K_{X}+L+ \frac{1}{\ell}A\mid_{S},h_{S}\cdot h_{A}^{\frac{1}{\ell}})^{n}$ (see Definition \ref{rmwp}), 
\begin{equation}\label{up2}
\limsup_{m\rightarrow\infty} \left((m!)^{-n}\cdot \prod_{j=1}^{m}d_{\ell,j}\right)^{\frac{1}{m}}\leqq \frac{1}{n!} 
(K_{X}+L+ \frac{1}{\ell}A\mid_{S},h_{S}\cdot h_{A}^{\frac{1}{\ell}})^{n}
\end{equation}
holds.
Hence combining (\ref{up1}) and (\ref{up2}), we see that 
\begin{eqnarray}\label{up3}
K_{\ell,\infty}& := & \limsup_{m\rightarrow\infty}\sqrt[n]{(m!)^{-n}K_{\ell,m}}
\\
& \leqq  &\frac{C}{n!}\cdot (K_{X}+L+ \frac{1}{\ell}A\mid_{S},h_{S}\cdot h_{A}^{\frac{1}{\ell}})^{n}\cdot dV^{-1}\cdot h_{L,0}^{-1}\cdot h_{A}^{-\frac{1}{\ell}} \nonumber
\end{eqnarray}
hold. 
\subsection{Lower estimate}
Now we shall estimate $K_{\ell,m}$ from below. 
Let $\nu$ denotes the numerical Kodaira dimension of 
$(K_{X}+L\mid_{S},h_{S})$. 
Let $V$ be a very general smooth complete intersection of $(n-\nu)$ members of 
 $\mid A\mid_{S}\mid$ such that \\$(K_{X}+L\mid_{S},h_{S})^{\nu}\cdot V > 0$. 

Then as in the proof of Lemma \ref{lower} in Section \ref{DyAZD},  we 
see that 
\begin{equation}\label{lo1}
\limsup_{\ell\rightarrow\infty}\ell^{n-\nu}\!\!\cdot K_{\ell,\infty} \neq 0 
\end{equation}
holds. 
Here the factor $\ell^{n-\nu}$ appears  by the decay of the curvature
of $\Theta_{h_{S}} + \ell^{-1}\Theta_{h_{A}}$  along $V$ in the normal direction of $V$.
Let $h_{V}$ be a singular hermitian metric on $K_{X} + L\mid_{V}$ 
such that 
\begin{enumerate}
\item $h_{V} \geqq h_{S}\!\!\mid_{V}$ holds on $V$. 
\item $\Theta_{h_{V}}$ is strictly positive everywhere on $V$. 
\end{enumerate}
Then by repeating  the proof of Lemma \ref{lower},
we see that  there exists a positive constant $C_{\varepsilon}$ depending only on $0 < \varepsilon < 1$ 
such that 
\begin{equation}\label{onV}
\limsup_{\ell\rightarrow\infty}\ell^{n-\nu}\!\!\cdot K_{\ell,\infty} \geqq 
C_{\varepsilon}\cdot (h_{V}^{\varepsilon}\cdot h_{S}\mid_{V}^{1-\varepsilon})^{-1}
\end{equation}
holds on $V$. 
\subsection{Completion of the proofs of Theorems \ref{main} and \ref{numsubad2} }

Let us set 
\[
K_{S,\infty}:= \limsup_{\ell\rightarrow\infty}
\left(\{(K_{X}+L+ \frac{1}{\ell}A,h_{S}\cdot h_{A}^{\frac{1}{\ell}})^{n}\}^{-1}\cdot (K_{\ell,\infty})
\right).
\]
Then since there exists a positive constant $c_{0}$ such that 
\[
(K_{X}+L+ \frac{1}{\ell}A,h_{S}\cdot h_{A}^{\frac{1}{\ell}})^{n}
\geqq c_{0}\cdot \left(\frac{1}{\ell^{n-\nu}}\right)
\]
holds by the assumption, 
by  (\ref{up3}) and (\ref{lo1}), we see that $K_{S,\infty}$ exists on $S$ and nonzero
and by (\ref{onV}), moving $V$, 
\[
h_{S,\infty} : = \mbox{the lower envelope of}\,\,\, \frac{1}{K_{S,\infty}}
\]
is an AZD of $(K_{X}+L\mid_{S},e^{-\varphi}\cdot dV^{-1}\cdot h_{L})$.

\noindent Repeating the argument in Section \ref{Dycon}, we may extend
$K_{\ell,\infty}$ to $\tilde{K}_{\ell,\infty}$ on $X$.
Then we set 
\[
\tilde{K}_{S,\infty}:= \limsup_{\ell\rightarrow\infty}
\left(\{(K_{X}+L+ \frac{1}{\ell}A\mid_{S},h_{S}\cdot h_{A}^{\frac{1}{\ell}})^{n}\}^{-1}\cdot (\tilde{K}_{\ell,\infty})
\right)
\]  
and 
\[
\tilde{h}_{S,\infty} := \mbox{the lower envelope of}\,\,\,\,\frac{1}{\tilde{K}_{S,\infty}}
\]
exists.  By the construction, it is clear that $\Theta_{\tilde{h}_{S,\infty}}$
is closed semipositive in the sense of current.  
Then we see that 
\[
\tilde{h}_{S,\infty}\mid_{S} \leqq h_{S,\infty}
\]
and $\tilde{h}_{S,\infty}\mid_{S}$ is an AZD of $(K_{X}+L\mid_{S},e^{-\varphi}\cdot dV^{-1}\cdot h_{L})$.

Then by Theorem \ref{extension}, we may extend every element of 
\[
H^{0}(S,{\cal O}_{S}(m(K_{X}+L))\otimes {\cal I}(e^{-\varphi}\cdot h_{L}\cdot h_{S,\infty}^{m-1}))
\]
to an element of 
\[
H^{0}(X,{\cal O}_{X}(m(K_{X}+L))\otimes {\cal I}(h_{L}\cdot \tilde{h}_{S,\infty}^{m-1})). 
\]
This completes the proof of Theorem \ref{main}. 
$\square$   \vspace{5mm} \\
{\bf Proof of Theorem \ref{numsubad2}.}
The proof of Theorem \ref{numsubad2} follows from the combination of 
the proof of Theorem \ref{numsubad} and the above argument of eliminating $A$.
$\square$

\small{

}
\noindent Author's address\\
Hajime Tsuji\\
Department of Mathematics\\
Sophia University\\
7-1 Kioicho, Chiyoda-ku 102-8554\\
Japan
\end{document}